\documentclass[8pt,a4wide]{article}

\usepackage[normalem]{ulem}
\usepackage{listings}
\usepackage{showlabels}
\usepackage{a4wide}
\usepackage{amsmath} 
\usepackage{amsthm}
\usepackage{amsfonts}
\usepackage{amssymb,amscd}
\usepackage{multirow}
\usepackage[mathcal]{eucal}
\usepackage[spanish,USenglish]{babel} 
\usepackage{graphics}  
\usepackage{graphicx} 
\usepackage{hyperref}
\usepackage{fancyhdr}
\usepackage{titlesec}
\usepackage[latin1]{inputenc}
\usepackage{amsmath}
\usepackage{amsfonts}
\usepackage{amssymb}
\usepackage{amsthm}
\usepackage{newlfont}
\usepackage{amsbsy}
\usepackage{bm}
\usepackage[all]{xy}
\usepackage{textcomp}
\usepackage{color}
\usepackage{epsfig}
\usepackage{dsfont}
\usepackage{latexsym}
\usepackage{array}
\usepackage{longtable}
\usepackage{enumerate}
\usepackage{multicol}
\usepackage{mathrsfs}
\usepackage{cancel}
\usepackage{wasysym}
\theoremstyle{plain}
\newtheorem{teor}{Theorem}
\newtheorem{lemma}{Lemma}
\newtheorem{coro}{Corollary}
\newtheorem{prop}{Proposition}

\theoremstyle{definition}

\newtheorem*{ejem*}{Examples}
\newtheorem{defin}{Definition}
\newtheorem*{demos}{Proof}
\newtheorem{remark}{Remark}
\theoremstyle{remark}

\newcommand{\complex}{\mathbf{\mathbb{C}}}

\newcommand{\natu}{\mathbf{\mathbb{N}}}
\newcommand{\llave}[1]{\left\{ #1\right\}}
\newcommand{\norma}[1]{\left \| #1 \right \|}
\newcommand{\vabs}[1]{\left| #1\right|}

\newcommand{\corch}[1]{\left[ #1\right]}
\newcommand{\paren}[1]{\left( #1\right)}
\newcommand{\QED}{\hfill \ensuremath{\Box}}

\usepackage[breakable]{tcolorbox}
\usepackage{parskip} 

\usepackage{iftex}
\ifPDFTeX
\usepackage[T1]{fontenc}
\usepackage{mathpazo}
\else
\usepackage{fontspec}
\fi

\usepackage{graphicx}

\usepackage{caption}
\DeclareCaptionFormat{nocaption}{}
\usepackage[Export]{adjustbox} 
\adjustboxset{max size={0.9\linewidth}{0.9\paperheight}}
\usepackage{float}
\usepackage{xcolor} 
\usepackage{enumerate} 
\usepackage{geometry} 
\usepackage{amsmath} 
\usepackage{amssymb} 
\usepackage{textcomp} 
\AtBeginDocument{%
}
\usepackage{upquote} 
\usepackage{eurosym} 
\usepackage[mathletters]{ucs} 
\usepackage{fancyvrb} 
\usepackage{grffile} 
\makeatletter 
\def\Gread@@xetex#1{%
	\IfFileExists{"\Gin@base".bb}%
	{\Gread@eps{\Gin@base.bb}}%
	{\Gread@@xetex@aux#1}%
}
\makeatother

\usepackage{hyperref}
\usepackage{titling}
\usepackage{longtable} 
\usepackage{booktabs}  
\usepackage[inline]{enumitem} 
\usepackage[normalem]{ulem} 
\usepackage{mathrsfs}

\definecolor{urlcolor}{rgb}{0,.145,.698}
\definecolor{linkcolor}{rgb}{.71,0.21,0.01}
\definecolor{citecolor}{rgb}{.12,.54,.11}

\definecolor{ansi-black}{HTML}{3E424D}
\definecolor{ansi-black-intense}{HTML}{282C36}
\definecolor{ansi-red}{HTML}{E75C58}
\definecolor{ansi-red-intense}{HTML}{B22B31}
\definecolor{ansi-green}{HTML}{00A250}
\definecolor{ansi-green-intense}{HTML}{007427}
\definecolor{ansi-yellow}{HTML}{DDB62B}
\definecolor{ansi-yellow-intense}{HTML}{B27D12}
\definecolor{ansi-blue}{HTML}{208FFB}
\definecolor{ansi-blue-intense}{HTML}{0065CA}
\definecolor{ansi-magenta}{HTML}{D160C4}
\definecolor{ansi-magenta-intense}{HTML}{A03196}
\definecolor{ansi-cyan}{HTML}{60C6C8}
\definecolor{ansi-cyan-intense}{HTML}{258F8F}
\definecolor{ansi-white}{HTML}{C5C1B4}
\definecolor{ansi-white-intense}{HTML}{A1A6B2}
\definecolor{ansi-default-inverse-fg}{HTML}{FFFFFF}
\definecolor{ansi-default-inverse-bg}{HTML}{000000}

\let\Oldtex\TeX
\let\Oldlatex\LaTeX
\renewcommand{\TeX}{\textrm{\Oldtex}}
\renewcommand{\LaTeX}{\textrm{\Oldlatex}}
\title{First Example}

\makeatletter
\def\PY@reset{\let\PY@it=\relax \let\PY@bf=\relax%
	\let\PY@ul=\relax \let\PY@tc=\relax%
	\let\PY@bc=\relax \let\PY@ff=\relax}
\def\PY@tok#1{\csname PY@tok@#1\endcsname}
\def\PY@toks#1+{\ifx\relax#1\empty\else%
	\PY@tok{#1}\expandafter\PY@toks\fi}
\def\PY@do#1{\PY@bc{\PY@tc{\PY@ul{%
				\PY@it{\PY@bf{\PY@ff{#1}}}}}}}
\def\PY#1#2{\PY@reset\PY@toks#1+\relax+\PY@do{#2}}

\expandafter\def\csname PY@tok@w\endcsname{\def\PY@tc##1{\textcolor[rgb]{0.73,0.73,0.73}{##1}}}
\expandafter\def\csname PY@tok@c\endcsname{\let\PY@it=\textit\def\PY@tc##1{\textcolor[rgb]{0.25,0.50,0.50}{##1}}}
\expandafter\def\csname PY@tok@cp\endcsname{\def\PY@tc##1{\textcolor[rgb]{0.74,0.48,0.00}{##1}}}
\expandafter\def\csname PY@tok@k\endcsname{\let\PY@bf=\textbf\def\PY@tc##1{\textcolor[rgb]{0.00,0.50,0.00}{##1}}}
\expandafter\def\csname PY@tok@kp\endcsname{\def\PY@tc##1{\textcolor[rgb]{0.00,0.50,0.00}{##1}}}
\expandafter\def\csname PY@tok@kt\endcsname{\def\PY@tc##1{\textcolor[rgb]{0.69,0.00,0.25}{##1}}}
\expandafter\def\csname PY@tok@o\endcsname{\def\PY@tc##1{\textcolor[rgb]{0.40,0.40,0.40}{##1}}}
\expandafter\def\csname PY@tok@ow\endcsname{\let\PY@bf=\textbf\def\PY@tc##1{\textcolor[rgb]{0.67,0.13,1.00}{##1}}}
\expandafter\def\csname PY@tok@nb\endcsname{\def\PY@tc##1{\textcolor[rgb]{0.00,0.50,0.00}{##1}}}
\expandafter\def\csname PY@tok@nf\endcsname{\def\PY@tc##1{\textcolor[rgb]{0.00,0.00,1.00}{##1}}}
\expandafter\def\csname PY@tok@nc\endcsname{\let\PY@bf=\textbf\def\PY@tc##1{\textcolor[rgb]{0.00,0.00,1.00}{##1}}}
\expandafter\def\csname PY@tok@nn\endcsname{\let\PY@bf=\textbf\def\PY@tc##1{\textcolor[rgb]{0.00,0.00,1.00}{##1}}}
\expandafter\def\csname PY@tok@ne\endcsname{\let\PY@bf=\textbf\def\PY@tc##1{\textcolor[rgb]{0.82,0.25,0.23}{##1}}}
\expandafter\def\csname PY@tok@nv\endcsname{\def\PY@tc##1{\textcolor[rgb]{0.10,0.09,0.49}{##1}}}
\expandafter\def\csname PY@tok@no\endcsname{\def\PY@tc##1{\textcolor[rgb]{0.53,0.00,0.00}{##1}}}
\expandafter\def\csname PY@tok@nl\endcsname{\def\PY@tc##1{\textcolor[rgb]{0.63,0.63,0.00}{##1}}}
\expandafter\def\csname PY@tok@ni\endcsname{\let\PY@bf=\textbf\def\PY@tc##1{\textcolor[rgb]{0.60,0.60,0.60}{##1}}}
\expandafter\def\csname PY@tok@na\endcsname{\def\PY@tc##1{\textcolor[rgb]{0.49,0.56,0.16}{##1}}}
\expandafter\def\csname PY@tok@nt\endcsname{\let\PY@bf=\textbf\def\PY@tc##1{\textcolor[rgb]{0.00,0.50,0.00}{##1}}}
\expandafter\def\csname PY@tok@nd\endcsname{\def\PY@tc##1{\textcolor[rgb]{0.67,0.13,1.00}{##1}}}
\expandafter\def\csname PY@tok@s\endcsname{\def\PY@tc##1{\textcolor[rgb]{0.73,0.13,0.13}{##1}}}
\expandafter\def\csname PY@tok@sd\endcsname{\let\PY@it=\textit\def\PY@tc##1{\textcolor[rgb]{0.73,0.13,0.13}{##1}}}
\expandafter\def\csname PY@tok@si\endcsname{\let\PY@bf=\textbf\def\PY@tc##1{\textcolor[rgb]{0.73,0.40,0.53}{##1}}}
\expandafter\def\csname PY@tok@se\endcsname{\let\PY@bf=\textbf\def\PY@tc##1{\textcolor[rgb]{0.73,0.40,0.13}{##1}}}
\expandafter\def\csname PY@tok@sr\endcsname{\def\PY@tc##1{\textcolor[rgb]{0.73,0.40,0.53}{##1}}}
\expandafter\def\csname PY@tok@ss\endcsname{\def\PY@tc##1{\textcolor[rgb]{0.10,0.09,0.49}{##1}}}
\expandafter\def\csname PY@tok@sx\endcsname{\def\PY@tc##1{\textcolor[rgb]{0.00,0.50,0.00}{##1}}}
\expandafter\def\csname PY@tok@m\endcsname{\def\PY@tc##1{\textcolor[rgb]{0.40,0.40,0.40}{##1}}}
\expandafter\def\csname PY@tok@gh\endcsname{\let\PY@bf=\textbf\def\PY@tc##1{\textcolor[rgb]{0.00,0.00,0.50}{##1}}}
\expandafter\def\csname PY@tok@gu\endcsname{\let\PY@bf=\textbf\def\PY@tc##1{\textcolor[rgb]{0.50,0.00,0.50}{##1}}}
\expandafter\def\csname PY@tok@gd\endcsname{\def\PY@tc##1{\textcolor[rgb]{0.63,0.00,0.00}{##1}}}
\expandafter\def\csname PY@tok@gi\endcsname{\def\PY@tc##1{\textcolor[rgb]{0.00,0.63,0.00}{##1}}}
\expandafter\def\csname PY@tok@gr\endcsname{\def\PY@tc##1{\textcolor[rgb]{1.00,0.00,0.00}{##1}}}
\expandafter\def\csname PY@tok@ge\endcsname{\let\PY@it=\textit}
\expandafter\def\csname PY@tok@gs\endcsname{\let\PY@bf=\textbf}
\expandafter\def\csname PY@tok@gp\endcsname{\let\PY@bf=\textbf\def\PY@tc##1{\textcolor[rgb]{0.00,0.00,0.50}{##1}}}
\expandafter\def\csname PY@tok@go\endcsname{\def\PY@tc##1{\textcolor[rgb]{0.53,0.53,0.53}{##1}}}
\expandafter\def\csname PY@tok@gt\endcsname{\def\PY@tc##1{\textcolor[rgb]{0.00,0.27,0.87}{##1}}}
\expandafter\def\csname PY@tok@err\endcsname{\def\PY@bc##1{\setlength{\fboxsep}{0pt}\fcolorbox[rgb]{1.00,0.00,0.00}{1,1,1}{\strut ##1}}}
\expandafter\def\csname PY@tok@kc\endcsname{\let\PY@bf=\textbf\def\PY@tc##1{\textcolor[rgb]{0.00,0.50,0.00}{##1}}}
\expandafter\def\csname PY@tok@kd\endcsname{\let\PY@bf=\textbf\def\PY@tc##1{\textcolor[rgb]{0.00,0.50,0.00}{##1}}}
\expandafter\def\csname PY@tok@kn\endcsname{\let\PY@bf=\textbf\def\PY@tc##1{\textcolor[rgb]{0.00,0.50,0.00}{##1}}}
\expandafter\def\csname PY@tok@kr\endcsname{\let\PY@bf=\textbf\def\PY@tc##1{\textcolor[rgb]{0.00,0.50,0.00}{##1}}}
\expandafter\def\csname PY@tok@bp\endcsname{\def\PY@tc##1{\textcolor[rgb]{0.00,0.50,0.00}{##1}}}
\expandafter\def\csname PY@tok@fm\endcsname{\def\PY@tc##1{\textcolor[rgb]{0.00,0.00,1.00}{##1}}}
\expandafter\def\csname PY@tok@vc\endcsname{\def\PY@tc##1{\textcolor[rgb]{0.10,0.09,0.49}{##1}}}
\expandafter\def\csname PY@tok@vg\endcsname{\def\PY@tc##1{\textcolor[rgb]{0.10,0.09,0.49}{##1}}}
\expandafter\def\csname PY@tok@vi\endcsname{\def\PY@tc##1{\textcolor[rgb]{0.10,0.09,0.49}{##1}}}
\expandafter\def\csname PY@tok@vm\endcsname{\def\PY@tc##1{\textcolor[rgb]{0.10,0.09,0.49}{##1}}}
\expandafter\def\csname PY@tok@sa\endcsname{\def\PY@tc##1{\textcolor[rgb]{0.73,0.13,0.13}{##1}}}
\expandafter\def\csname PY@tok@sb\endcsname{\def\PY@tc##1{\textcolor[rgb]{0.73,0.13,0.13}{##1}}}
\expandafter\def\csname PY@tok@sc\endcsname{\def\PY@tc##1{\textcolor[rgb]{0.73,0.13,0.13}{##1}}}
\expandafter\def\csname PY@tok@dl\endcsname{\def\PY@tc##1{\textcolor[rgb]{0.73,0.13,0.13}{##1}}}
\expandafter\def\csname PY@tok@s2\endcsname{\def\PY@tc##1{\textcolor[rgb]{0.73,0.13,0.13}{##1}}}
\expandafter\def\csname PY@tok@sh\endcsname{\def\PY@tc##1{\textcolor[rgb]{0.73,0.13,0.13}{##1}}}
\expandafter\def\csname PY@tok@s1\endcsname{\def\PY@tc##1{\textcolor[rgb]{0.73,0.13,0.13}{##1}}}
\expandafter\def\csname PY@tok@mb\endcsname{\def\PY@tc##1{\textcolor[rgb]{0.40,0.40,0.40}{##1}}}
\expandafter\def\csname PY@tok@mf\endcsname{\def\PY@tc##1{\textcolor[rgb]{0.40,0.40,0.40}{##1}}}
\expandafter\def\csname PY@tok@mh\endcsname{\def\PY@tc##1{\textcolor[rgb]{0.40,0.40,0.40}{##1}}}
\expandafter\def\csname PY@tok@mi\endcsname{\def\PY@tc##1{\textcolor[rgb]{0.40,0.40,0.40}{##1}}}
\expandafter\def\csname PY@tok@il\endcsname{\def\PY@tc##1{\textcolor[rgb]{0.40,0.40,0.40}{##1}}}
\expandafter\def\csname PY@tok@mo\endcsname{\def\PY@tc##1{\textcolor[rgb]{0.40,0.40,0.40}{##1}}}
\expandafter\def\csname PY@tok@ch\endcsname{\let\PY@it=\textit\def\PY@tc##1{\textcolor[rgb]{0.25,0.50,0.50}{##1}}}
\expandafter\def\csname PY@tok@cm\endcsname{\let\PY@it=\textit\def\PY@tc##1{\textcolor[rgb]{0.25,0.50,0.50}{##1}}}
\expandafter\def\csname PY@tok@cpf\endcsname{\let\PY@it=\textit\def\PY@tc##1{\textcolor[rgb]{0.25,0.50,0.50}{##1}}}
\expandafter\def\csname PY@tok@c1\endcsname{\let\PY@it=\textit\def\PY@tc##1{\textcolor[rgb]{0.25,0.50,0.50}{##1}}}
\expandafter\def\csname PY@tok@cs\endcsname{\let\PY@it=\textit\def\PY@tc##1{\textcolor[rgb]{0.25,0.50,0.50}{##1}}}

\makeatother

\makeatletter
\newbox\Wrappedcontinuationbox 
\newbox\Wrappedvisiblespacebox 
\newcommand*\Wrappedvisiblespace {\textcolor{red}{\textvisiblespace}} 
\newcommand*\Wrappedcontinuationsymbol {\textcolor{red}{\llap{\tiny$\m@th\hookrightarrow$}}} 
\newcommand*\Wrappedcontinuationindent {3ex } 
\newcommand*\Wrappedafterbreak {\kern\Wrappedcontinuationindent\copy\Wrappedcontinuationbox} 
\newcommand*\Wrappedbreaksatspecials {%
	\def\PYGZus{\discretionary{\char`\_}{\Wrappedafterbreak}{\char`\_}}%
	\def\PYGZob{\discretionary{}{\Wrappedafterbreak\char`\{}{\char`\{}}%
	\def\PYGZcb{\discretionary{\char`\}}{\Wrappedafterbreak}{\char`\}}}%
	\def\PYGZca{\discretionary{\char`\^}{\Wrappedafterbreak}{\char`\^}}%
	\def\PYGZam{\discretionary{\char`\&}{\Wrappedafterbreak}{\char`\&}}%
	\def\PYGZlt{\discretionary{}{\Wrappedafterbreak\char`\<}{\char`\<}}%
	\def\PYGZgt{\discretionary{\char`\>}{\Wrappedafterbreak}{\char`\>}}%
	\def\PYGZsh{\discretionary{}{\Wrappedafterbreak\char`\#}{\char`\#}}%
	\def\PYGZpc{\discretionary{}{\Wrappedafterbreak\char`\%}{\char`\%}}%
	\def\PYGZdl{\discretionary{}{\Wrappedafterbreak\char`\$}{\char`\$}}%
	\def\PYGZhy{\discretionary{\char`\-}{\Wrappedafterbreak}{\char`\-}}%
	\def\PYGZsq{\discretionary{}{\Wrappedafterbreak\textquotesingle}{\textquotesingle}}%
	\def\PYGZdq{\discretionary{}{\Wrappedafterbreak\char`\"}{\char`\"}}%
	\def\PYGZti{\discretionary{\char`\~}{\Wrappedafterbreak}{\char`\~}}%
} 
\newcommand*\Wrappedbreaksatpunct {%
	\lccode`\~`\.\lowercase{\def~}{\discretionary{\hbox{\char`\.}}{\Wrappedafterbreak}{\hbox{\char`\.}}}%
	\lccode`\~`\,\lowercase{\def~}{\discretionary{\hbox{\char`\,}}{\Wrappedafterbreak}{\hbox{\char`\,}}}%
	\lccode`\~`\;\lowercase{\def~}{\discretionary{\hbox{\char`\;}}{\Wrappedafterbreak}{\hbox{\char`\;}}}%
	\lccode`\~`\:\lowercase{\def~}{\discretionary{\hbox{\char`\:}}{\Wrappedafterbreak}{\hbox{\char`\:}}}%
	\lccode`\~`\?\lowercase{\def~}{\discretionary{\hbox{\char`\?}}{\Wrappedafterbreak}{\hbox{\char`\?}}}%
	\lccode`\~`\!\lowercase{\def~}{\discretionary{\hbox{\char`\!}}{\Wrappedafterbreak}{\hbox{\char`\!}}}%
	\lccode`\~`\/\lowercase{\def~}{\discretionary{\hbox{\char`\/}}{\Wrappedafterbreak}{\hbox{\char`\/}}}%
	\catcode`\.\active
	\catcode`\,\active 
	\catcode`\;\active
	\catcode`\:\active
	\catcode`\?\active
	\catcode`\!\active
	\catcode`\/\active 
	\lccode`\~`\~ 	
}
\makeatother

\let\OriginalVerbatim=\Verbatim
\makeatletter
\renewcommand{\Verbatim}[1][1]{%
	\sbox\Wrappedcontinuationbox {\Wrappedcontinuationsymbol}%
	\sbox\Wrappedvisiblespacebox {\FV@SetupFont\Wrappedvisiblespace}%
	\def\FancyVerbFormatLine ##1{\hsize\linewidth
		\vtop{\raggedright\hyphenpenalty\z@\exhyphenpenalty\z@
			\doublehyphendemerits\z@\finalhyphendemerits\z@
			\strut ##1\strut}%
	}%
	\def\FV@Space {%
		\nobreak\hskip\z@ plus\fontdimen3\font minus\fontdimen4\font
		\discretionary{\copy\Wrappedvisiblespacebox}{\Wrappedafterbreak}
		{\kern\fontdimen2\font}%
	}%
	
	\Wrappedbreaksatspecials
	\OriginalVerbatim[#1,codes*=\Wrappedbreaksatpunct]%
}
\makeatother

\definecolor{incolor}{HTML}{303F9F}
\definecolor{outcolor}{HTML}{D84315}
\definecolor{cellborder}{HTML}{CFCFCF}
\definecolor{cellbackground}{HTML}{F7F7F7}

\makeatletter
\newcommand{\boxspacing}{\kern\kvtcb@left@rule\kern\kvtcb@boxsep}
\makeatother

\hypersetup{
	breaklinks=true,  
	colorlinks=true,
	urlcolor=urlcolor,
	linkcolor=linkcolor,
	citecolor=citecolor,
}

\begin{document}
	\title{\textbf{A Class of Matrix Schr\"{o}dinger Bispectral Operators}}

 	\author{
Brian D. Vasquez Campos$^{1}$\\ \\
\small{$^{1}$Department of Mathematics, Khalifa University.}\\
\small{\texttt{$^{1}$bridava927@gmail.com}}
}

\maketitle


\begin{abstract}
	We prove the bispectrality of some class of matrix Schr\"{o}dinger operators with polynomial potentials which satisfy a second-order matrix autonomous differential equation. The physical equation is constructed using the formal theory of the Laurent series and after that obtaining local solutions using estimations in the Frobenius norm. Furthermore, the characterization of the 
	algebra of polynomial eigenvalues in the spectral variable is given using some family of functions $\mathcal{P}=\llave{P_{k}}_{k\in \natu}$  with the remarkable property of satisfying a general version of the Leibniz rule.

\end{abstract}

\section*{Introduction}

The Bispectral Problem was originally posed by J.~J.~Duistermaat and F.~A.~Gr\"{u}nbaum \cite{duistermaat1986differential}, it consists in to find all the bispectral triples $(L,\psi, B)$ 
satisfying systems of equations 
\begin{equation}\label{bispectral}
L\psi(x,z)=\psi (x,z)F(z) \hspace{1 cm}  (\psi B) (x,z)=\theta(x)\psi(x,z)
\end{equation}
with $L=L(x,\partial_{x})$, $B=B(z,\partial_{z})$  linear matrix  differential operators, i.e., $L\psi=\sum_{i=0}^{l}a_{i}(x) \cdot \partial_{x}^i \psi$,
$\psi B=\sum_{j=0}^{m}\partial_{z}^j \psi \cdot b_{j}(z)$. The functions $a_{i}, b_{j},F, \theta$ and the nontrivial common eigenfunction $\psi$ are in principle compatible sized matrix valued functions.

The bispectral problem was completely solved in the scalar case for Schr\"{o}dinger  operators $L=-\partial_{x}^2+U(x)$ and the potentials $U(x)$ for which bispectrality follows were characterized in
 \cite{duistermaat1986differential}. It turned out to have deep connections with many problems in Mathematical Physics.  Indeed, it could be 
 arranged in suitable manifolds which were naturally parameterized by the flows of the Korteweg de-Vries (KdV) hierarchy or its master-symmetries
 \cite{zubelli1991differential,Zubelli2000}. It led to generalizations associated to the  Kadomtsev-Petviashvili (KP) hierarchy~\cite{Zubelli1992,wilson,horozov2002bispectral,iliev1999discrete}.

In the references \cite{VZ1, VZ2, GVZ}, we elucidate the symmetry structure inherent in a non-commutative variant of the bispectral problem. Additionally, we present a comprehensive framework and establish a method for determining the bispectrality of a given matrix polynomial sub-algebra within this context.

In this paper we construct the physical equation $L\psi= -z^2 \psi$ with the  matrix Schr\"{o}dinger  operator \\ $L=-\partial_{x}^2+V'(x)$ and the eigenfunction 
\begin{equation*}
\psi (x,z)=\paren{Iz+\frac{1}{2}V(x)}e^{xz}.
\end{equation*}

The necessary and sufficient condition to obtain the physical equation in this particular case is \\ $V''(x)=V'(x)V(x)$.

The main goals of this article are: Firstly, to obtain meromorphic solutions of the physical equation using the theory of the Laurent series. Secondly, we give a characterization of the algebra of polynomial eigenvalues in the spectral variable independently of the bispectrality and look for conditions on the function $V$  to obtain that this algebra is not trivial, this is precisely bispectrality.

The bispectral problem over a noncommutative ring was first studied in \cite{zubelli1990differential, Zubelli1992a, Zubelli1992b} in the case of the ring of matrices for the situation where both the physical and spectral operators were acting on the same side of the eigenfunction. Later on, several generalizations were considered. 
See~\cite{SZ01, Kas15, BL08, GI03, Grfrm-4, GHY17} and references therein. This article follows up on the possibility of having the physical and the spectral operators acting on different sides. If the operator in the physical variable is normalized, its corresponding eigenvalue is scalar and the operators act from different sides we obtain the ad-condition and the right bispectral algebra is contained in the algebra of matrix polynomials $M_{N}(\complex[x])$.

The plan of this article is as follows: In Section \ref{auto}, we study this matrix autonomous differential equation using Laurent series with a simple pole at the origin $V(x)=\sum_{k=-1}^{\infty}V_{k}x^{k}$ and obtaining conditions on the coefficients which gives rise to some remarkable properties such as  $V_{k}(V_{0},V_{1},V_{2})$ is quasihomogeneous of type $(1,2,3)$ and degree $k+1$. After that, we obtain important estimations in the Frobenius norm to assure the existence of local meromorphic solutions of the autonomous equation $V''(x)=V'(x)V(x)$. An important property of working in the matrix case is the existence of nonconstant polynomial solutions of this autonomous equation. In section \ref{bisch}, we give a complete characterization of the algebra 
\begin{equation*}
\mathbb{A}=\llave{\theta \in M_{N}(\complex[x])\mid \exists B=B(z,\partial_{z}), (\psi B)(x,z) =\theta(x) \psi(x,z)}
\end{equation*}
using the family of functions $\mathcal{P}=\llave{P_{k}}_{k\in \natu}$ defined by 
\begin{equation*}
P_{k}(\theta)=\frac{\theta^{(k)}(0)}{(k-1)!}-\frac{1}{2}\sum_{j=0}^{k}\corch{\frac{\theta^{(j)}(0)}{j!},V_{k-1-j}},
\end{equation*}
where $k\in \natu$, $\theta \in M_{N}(\complex[x])$. Furthermore, we prove the bispectral property for some class of polynomial potentials satisfying $V''(x)=V'(x)V(x)$.

\section{Algebraic Morphisms Arising from the  Matrix Equation \texorpdfstring{$V^{\prime \prime}(x)=V^{\prime}(x)V(x)$}{Lg}}\label{auto}	

We try to find meromorphic solutions for the matrix equation  $V''(x)=V'(x)V(x)$ with a simple pole at $x=0$. This will allow us to obtain some type of solutions of the bispectral problem. Leading to obtain formal solutions of the autonomous equation we see that the Taylor coefficients in the expansion $V(x)=\sum_{k=-1}^{\infty}V_{k}x^{k}$ turns out to be affine algebraic morphims. Furthermore, if the residue $V_{-1}=Res(V,0)=0$, then the holomorphic solution $V$ has Taylor coefficients which are quasihomogeneous in the noncommutative variables $V_0 , V_1$ and using some grading we obtain the bispectrality in the case of polynomial potentials.

\subsection{The Matrix Equation $V^{\prime \prime}(x)=V^{\prime}(x)V(x)$. }

Let  $V(x)=\sum_{k=-1}^{\infty}V_{k}x^{k}$, then $V^{'}(x)=\sum_{k=-1}^{\infty}kV_{k}x^{k-1}$ and $V^{\prime \prime}(x)=\sum_{k=-1}^{\infty}k(k-1)V_{k}x^{k-2}$. Therefore, $V^{\prime \prime}(x)=V^{\prime}(x)V(x)$ if, and only if, 
\begin{equation*}
k(k-1)V_k = \sum_{j=-1}^{k}jV_{j}V_{k-1-j}
\end{equation*}
for $k=-1,0,1,\cdots$.
\begin{itemize}
	\item If $k=-1$, then $-V_{-1}^2=2V_{-1}$ and hence $V_{-1}(V_{-1}+2I_{N})=0$.
	
	\item 	If $k=0$, then $V_{-1}V_{0}=0$.
	
	\item 	If $k=1$, then $V_{-1}V_{1}=0$.
	
	\item 	If $k\geq 2$, then $k(k-1)V_k = -V_{-1}V_k +kV_{k}V_{-1}+ \sum_{j=1}^{k-1} j V_{j}V_{k-1-j}$. Thus, 
	\begin{equation*}
	T_{k}(V_{k})= \sum_{j=1}^{k-1} j V_{j}V_{k-1-j}.
	\end{equation*}
where the operator $T_{k}:M_{N}(\complex)\rightarrow M_{N}(\complex)$ defined by $T_{k}(a)=k(k-1)a+V_{-1}a-ka V_{-1}$ \mbox{ .}
\end{itemize}
Since $V_{-1}(V_{-1}+2I_N)=0$, we have that $0$ and $-2$ are the only eigenvalues of $V_{-1}$. The Jordan Canonical Form Theorem implies that 
$V_{-1}$ has the form $diag(-2,0,\cdots , 0, -2)$. After a change of coordinates, we may  assume without loss of generality that
\begin{equation*}
V_{-1}=\paren{
	\begin{matrix}  
	-2I_{m} & 0 \\
	0 & 0
	\end{matrix}}.
\end{equation*}
Since $V_{-1}V_{0}=0$ we have that 

\begin{equation*}
\paren{
	\begin{matrix}  
	0  & 0 \\
	0  & 0
	\end{matrix}}
=\paren{
	\begin{matrix}  
	-2I_{m} & 0 \\
	0 & 0
	\end{matrix}}
\paren{
	\begin{matrix}  
	V_{011}  & V_{012} \\
	V_{021}  & V_{022}
	\end{matrix}}
=\paren{
	\begin{matrix}  
	-2V_{011}  & -2V_{012} \\
	0 & 0
	\end{matrix}}.
\end{equation*}
Then,  $V_{011}=V_{012}=0$. Thus, 
\begin{equation*}
V_{0}=\paren{
	\begin{matrix}  
	0 & 0 \\
	V_{021}  & V_{022}
	\end{matrix}}.
\end{equation*}
In the same way,   $V_{-1}V_{1}=0$ implies that

\begin{equation*}
V_{1}=\paren{
	\begin{matrix}  
	0 & 0 \\
	V_{121}  & V_{122}
	\end{matrix}}.
\end{equation*}
Now we write 
\begin{equation*}
V_{k}=\paren{
	\begin{matrix}  
	V_{k11}  & V_{k12} \\
	V_{k21}  & V_{k22}
	\end{matrix}},
\end{equation*}
to obtain 

\begin{equation*}
T_{k}(V_{k})=k(k-1)V_{k}+V_{-1}V_{k}-kV_{k} V_{-1}
=k(k-1)\paren{
	\begin{matrix}  
	V_{k11}  & V_{k12} \\
	V_{k21}  & V_{k22}
	\end{matrix}}+
\paren{
	\begin{matrix}  
	-2I_{m} & 0 \\
	0 & 0
	\end{matrix}}\paren{
	\begin{matrix}  
	V_{k11}  & V_{k12} \\
	V_{k21}  & V_{k22}
	\end{matrix}}
\end{equation*}
\begin{equation*}
-k\paren{
	\begin{matrix}  
	V_{k11}  & V_{k12} \\
	V_{k21}  & V_{k22}
	\end{matrix}}\paren{
	\begin{matrix}  
	-2I_{m} & 0 \\
	0 & 0
	\end{matrix}}
=k(k-1)\paren{
	\begin{matrix}  
	V_{k11}  & V_{k12} \\
	V_{k21}  & V_{k22}
	\end{matrix}}
+\paren{
	\begin{matrix}  
	-2V_{k11}  &-2 V_{k12} \\
	0  & 0
	\end{matrix}}
+k\paren{
	\begin{matrix}  
	2V_{k11}  & 0\\
	2V_{k21}  & 0
	\end{matrix}}
\end{equation*}
\begin{equation*}
=\paren{
	\begin{matrix}  
	(k(k-1)+2k-2)V_{k11}  & (k(k-1)-2)V_{k12} \\
	(k(k-1)+2k)V_{k21}  & k(k-1)V_{k22}
	\end{matrix}}
\end{equation*}
\begin{equation*}
=\paren{
	\begin{matrix}  
	(k-1)(k+2)V_{k11}  & (k-2)(k+1)V_{k12} \\
	k(k+1)V_{k21}  & k(k-1)V_{k22}
	\end{matrix}} \mbox{ .}
\end{equation*}
For $k=2$, we have 
\begin{equation*}
V_{1}V_{0}=
\paren{
	\begin{matrix}  
	0 & 0 \\
	V_{121}  & V_{122}
	\end{matrix}}
\paren{
	\begin{matrix}  
	0 & 0 \\
	V_{021}  & V_{022}
	\end{matrix}}
=\paren{
	\begin{matrix}  
	0 & 0 \\
	V_{122}V_{021}  & V_{122}V_{022}
	\end{matrix}}
= T_{2}(V_{2})=
\paren{
	\begin{matrix}  
	4V_{211}  & 0 \\
	6V_{221}  & 2V_{222}
	\end{matrix}}.
\end{equation*}
Therefore, $V_{211}=0$, $V_{221}=\frac{V_{122}V_{021}}{6}$, $V_{222}=\frac{V_{122}V_{022}}{2}$. Thus,

\begin{equation*}
V_{2}=\paren{
	\begin{matrix}  
	0  & V_{212} \\
	\frac{V_{122}V_{021}}{6} & \frac{V_{122}V_{022}}{2}
	\end{matrix}}
\end{equation*}
Remember that,
\begin{equation*}
T_{k}(V_{k})=k(k-1)V_{k}+V_{-1}V_{k}-kV_{k} V_{-1}
\end{equation*}
for $k\geq 2$. 
If $k\geq 3$, we have that $T_{k}$ is invertible and 
\begin{equation}\label{poly}
V_{k}= \sum_{j=1}^{k-1} j  T_{k}^{-1}(V_{j}V_{k-1-j}).
\end{equation}

\begin{defin}
	Fix an element $A\in M_{N}(\complex)$, define the multiplication operators $L_{A},R_{A}:M_{N}(\complex)\rightarrow M_{N}(\complex)$,
	$L_{A}(X)=AX$, $R_{A}(X)=XA$.
\end{defin}

We now look for the elements $A\in M_{N}(\complex)$ such that $L_{A}$ and $R_{A}$ commutes with $T_{k}$ for $k\geq 1$.

\begin{lemma}\label{conmu}
$L_{A}$ and $R_{A}$ commutes with $T_{k}$ for $k\geq 1$ if, and only if, $A_{12}=0\in M_{m\times (N-m)}(\complex)$ and $A_{21}=0\in M_{(N-m)\times m}(\complex)$.
\end{lemma}
\begin{demos}
	Note that 
	\begin{equation*}
	T_{k}R_{A}(X)= T_{k}(XA)=k(k-1)XA+V_{-1}XA-k(XA)V_{-1}
	=k(k-1)XA+V_{-1}XA-kXV_{-1}-kX[A,V_{-1}]
	\end{equation*}
		\begin{equation*}
	=T_{k}(X)A-kX[A,V_{-1}]=R_{A}T_{k}(X)-kX[A,V_{-1}],
		T_{k}L_{A}(X)=T_{k}(AX)=k(k-1)AX+V_{-1}AX-k(AX)V_{-1}
			\end{equation*}
		\begin{equation*}
		=k(k-1)AX+AV_{-1}X-k(AX)V_{-1}-[A,V_{-1}]X=AT_{k}(X)-[A,V_{-1}]X.
	\end{equation*} 
	 Then,  $L_{A}$ and $R_{A}$ commutes with $T_{k}$ for $k\geq 1$ if, and only if, 
	$[A,V_{-1}]=0$, but this condition says that $A_{12}=0\in M_{m\times (N-m)}(\complex)$ and $A_{21}=0\in M_{(N-m)\times m}(\complex)$. \QED
\end{demos}

The following lemma gives us an interesting property of the operator $T_{k}$ when we consider $M_{N}(\complex)$ with the Frobenius norm.

\begin{lemma}
	The operator $T_{k}^{-1}:M_{N}(\complex) \rightarrow M_{N}(\complex)$ satisfies 
	\begin{equation*}
	\norma{T_{k}^{-1}(a)}_{F}\leq \frac{4(k^2-3)}{(k-2)(k-1)(k+1)(k+2)}\norma{a}_{F}
	\end{equation*}
	for $k\geq 3$.
\end{lemma}
\begin{demos}
	Since, 
	\begin{equation*}
T_{k}(a)	=\paren{
		\begin{matrix}  
		(k-1)(k+2)a_{11}  & (k-2)(k+1)a_{12} \\
		k(k+1)a_{21}  & k(k-1)a_{22}
		\end{matrix}}
	\end{equation*}
	for $k\geq 3$.
We have,
	\begin{equation*}
	T_{k}^{-1}(a)	=\paren{
	\begin{matrix}  
\frac{1}{(k-1)(k+2)} a_{11}  & \frac{1}{(k-2)(k+1)} a_{12} \\
\frac{1}{k(k+1)} a_{21}  & \frac{1}{k(k-1)} a_{22}
	\end{matrix}}.
\end{equation*}
Applying the Frobenius norm, 
\begin{equation*}
\norma{T_{k}^{-1}(a)}_{F}^2 =
\paren{\frac{1}{(k-1)(k+2)}}^2\norma{a_{11}}_{F}^2 +
\paren{\frac{1}{(k-2)(k+1)}}^2\norma{a_{12}}_{F}^2 +
\paren{\frac{1}{k(k+1)}}^2\norma{a_{21}}_{F}^2
\end{equation*}
\begin{equation*}
+
\paren{\frac{1}{k(k-1)}}^2\norma{a_{22}}_{F}^2 
\leq \llave{ \paren{\frac{1}{(k-1)(k+2)}}^2+\paren{\frac{1}{(k-2)(k+1)}}^2+\paren{\frac{1}{k(k+1)}}^2
+\paren{\frac{1}{k(k-1)}}^2}\norma{a}_{F}^2.
\end{equation*}
Therefore,
\begin{equation*}
\norma{T_{k}^{-1}(a)}_{F}\leq 
\norma{\paren{\frac{1}{(k-1)(k+2)},\frac{1}{(k-2)(k+1)},\frac{1}{k(k+1)},\frac{1}{k(k-1)}}}_{2}
\norma{a}_{F}
\end{equation*}
\begin{equation*}
= 2 \frac{\sqrt{k^6-5k^4+6k^2+8}}{(k-2)(k-1)k(k+1)(k+2)}\norma{a}_{F}
\leq \frac{4(k^2-3)}{(k-2)(k-1)(k+1)(k+2)}\norma{a}_{F},
\end{equation*}
for $k\geq 3$.
To obtain the last inequality observe that
\begin{equation*}
    \frac{\sqrt{k^6-5k^4+6k^2+8}}{(k-2)(k-1)k(k+1)(k+2)}\leq 
    \frac{2(k^2-3)}{(k-2)(k-1)(k+1)(k+2)}
\end{equation*}
\begin{equation*}
\Longleftrightarrow \sqrt{k^6-5k^4+6k^2+8}\leq 2k(k^2-3)
\end{equation*}
\begin{equation*}
\Longleftrightarrow k^6-5k^4+6k^2+8\leq 4k^2 (k^2-3)^2 
=4k^2(k^4-6k^2+9) 
=4k^6-24k^4+36k^2
\end{equation*}
\begin{equation*}
\Longleftrightarrow f(k):=3k^6-19k^4+30k^2-8\geq 0.
\end{equation*}
However, we have the factorization of the polynomial $f\in \mathbb{Q}[x]$, $f(x)=(x^2-4)(x^2-2)(3x^2-1)$, and $f(x)\geq 0$ for $\vabs{x}\geq 2$.\QED
\end{demos}

\begin{remark}
Note that the inequality in the previous lemma implies that $T_{k}^{-1}$ is a contraction for $k\geq 3$.
\end{remark}

We can use this result to estimate the norm of the sequence $\llave{V_{j}}_{j\in \natu}$.

\begin{teor}\label{convergence}
	If $\norma{V_{0}}_{F}\leq \frac{1}{4}$, $\norma{V_{1}}_{F}\leq \frac{1}{8}$, $\norma{V_{2}}_{F}\leq \frac{1}{16}$, then 
	\begin{equation*}
	\norma{V_{k}}_{F}\leq \frac{1}{2^{k+2}}
	\end{equation*}
	for every $k\geq 3$.
	\end{teor}

\begin{demos}
The proof is by induction. By hypothesis we have  $\norma{V_{k}}_{F}\leq \frac{1}{2^{k+2}}$ for $k=0,1,2$. Assume the claim for some $0\leq j\leq k-1$ and note that
\begin{equation*}
\norma{V_{k}}_{F}\leq \sum_{j=1}^{k-1} j  \norma{T_{k}^{-1}(V_{j}V_{k-1-j})}_{F}
\leq \sum_{j=1}^{k-1} j  \frac{4(k^2-3)}{(k-2)(k-1)(k+1)(k+2)}\norma{V_{j}V_{k-1-j}}_{F}
\end{equation*}
\begin{equation*}
\leq \sum_{j=1}^{k-1} j  \frac{4(k^2-3)}{(k-2)(k-1)(k+1)(k+2)}\norma{V_{j}}_{F}\norma{V_{k-1-j}}_{F}
\end{equation*}
\begin{equation*}
\leq \sum_{j=1}^{k-1} j  \frac{4(k^2-3)}{(k-2)(k-1)(k+1)(k+2)}\paren{\frac{1}{2^{j+2}}}\paren{\frac{1}{2^{k-j+1}}}
\end{equation*}
\begin{equation*}
= \frac{k(k^2-3)}{(k-2)(k+1)(k+2)}\paren{\frac{1}{2^{k+2}}}\leq \frac{1}{2^{k+2}}.
\end{equation*}
Therefore, the claim follows by induction.\QED
\end{demos}

This theorem allows us to give meromorphic solutions to the matrix equation $V''(x)=V'(x)V(x)$ in a punctured neighborhood of the origin. To do this, we consider the set $K\subset 
M_{N}(\complex)^3  \times \complex$ defined by the relations 
\begin{equation*}
V_{-1}V_{0}=V_{-1}V_{1}=0,
	V_{1}V_{0}=T_{2}(V_{2}), \norma{V_{0}}_{F}\leq \frac{1}{4}, \norma{V_{1}}_{F}\leq \frac{1}{8}, \norma{V_{2}}_{F}\leq \frac{1}{16}, 0<\vabs{x}\leq 1 .
\end{equation*}

\begin{coro}
The formal power series $V=V(V_{0},V_{1},V_{2},x)=\sum_{k=-1}^{\infty}V_{k}(V_{0},V_{1},V_{2})x^{k}$ is meromorphic for $(V_{0},V_{1},V_{2},x)\in K$.
\end{coro}
\begin{demos}
	If $(V_{0},V_{1},V_{2},x)\in K$, then Theorem~\ref{convergence} implies that $\norma{V_{k}(V_{0},V_{1},V_{2})}_{F}\leq \frac{1}{2^{k+2}}$ and therefore 
	$\norma{V_{k}(V_{0},V_{1},V_{2})x^{k}}_{F}\leq \frac{1}{2^{k+2}}$. Thus, the Weiertrass Theorem implies that the series $\sum_{k=-1}^{\infty}V_{k}(V_{0},V_{1},V_{2})x^{k}$ converges 	absolutely and uniformly in compact subsets of $K$. Since the functions $V_{k}(V_{0},V_{1},V_{2})x^{k}$ are meromorphic in $K$ we obtain the same for $V$. \QED
	
\end{demos}
\subsection{Some Properties of the Sequence $\llave{V_{k}(V_{0},V_{1},V_{2})}_{k\in \natu}$ }

The sequence $\llave{V_{k}(V_{0},V_{1},V_{2})}_{k\in \natu}$ has  important properties which are given in the following results.
\begin{prop}
	The function  $V_{k}(V_{021},V_{022},V_{121},V_{122},V_{212})$ has polynomial coordinates for every $k\in \natu$.
\end{prop}
\begin{demos}
Note that 
\begin{equation*}
V_{0}=\paren{
	\begin{matrix}  
	0  & 0 \\
	V_{021} & V_{022}
	\end{matrix}}
\end{equation*}
\begin{equation*}
V_{1}=\paren{
	\begin{matrix}  
	0  & 0 \\
	V_{121} & V_{122}
	\end{matrix}}
\end{equation*}
and \begin{equation*}
V_{2}=\paren{
	\begin{matrix}  
	0  & V_{212} \\
	\frac{V_{122}V_{021}}{6} & \frac{V_{122}V_{022}}{2}
	\end{matrix}}
\end{equation*}
has polynomial coordinates in $V_{021},V_{022},V_{121},V_{122},V_{212}$. Assume that $V_{j}(V_{021},V_{022},V_{121},V_{122},V_{212})$ has polynomial coordinates 
for $1\leq j \leq k-1$, since 
\begin{equation*}
V_{k}=\sum_{j=1}^{k-1}jT_{k}^{-1}(V_{j}V_{k-1-j})
\end{equation*}
and the product of matrices $V_{j}V_{k-1-j}$ is a polynomial in the block entries of $V_{j}$ and $V_{k-1-j}$. Since $T_{k}^{-1}$ is linear we obtain that $V_{k}(V_{021},V_{022},V_{121},V_{122},V_{212})$ has polynomial coordinates. \QED

\end{demos}
\begin{coro}
	$V_{k}(V_{021},V_{022},V_{121},V_{122},V_{212})$ is an algebraic morphism for every $k\in \natu$.
\end{coro}

\begin{teor}
	If $A\in M_{N}(\complex)$ and $A_{22}V_{j,21}=V_{j,21}A_{11}$, for $j=0,1,2$, $A_{11}V_{12}=V_{12}A_{22}$ , $[A_{22},V_{0,22}]=[A_{22},V_{1,22}]=0$, then 
	\begin{equation*}
	V_{j}(V_{j21}A_{11},V_{j22}A_{22})=A^{j+1}V_{j}(V_{j21},V_{j22})=V_{j}(V_{j21},V_{j22}) A^{j+1},
	\end{equation*}
	for $j=0,1$.
	\begin{equation*}
	V_{2}(AV_{0},A^2 V_{1}, A_{11}^3 V_{212})=A^3V_{2}(V_{0},V_{1},V_{212})=V_{2}(V_{0},V_{1},V_{212})A^{3}
	\end{equation*}
	\begin{equation*}
		V_{k}(A V_{0},A^2 V_{1}, A^3 V_{2})=A^{k+1} V_{k}(V_{0},V_{1},V_{2})=V_{k}(V_{0},V_{1},V_{2})A^{k+1},
	\end{equation*}
	for every $k \geq 3$.
\end{teor}

\begin{demos}
In fact, 
\begin{equation*}
AV_{0}(V_{021},V_{022})=
\paren{
	\begin{matrix}  
	A_{11}  & 0 \\
	0 & A_{22}
	\end{matrix}}
\paren{
	\begin{matrix}  
	0 & 0 \\
	V_{021}  & V_{022}
	\end{matrix}}
=\paren{
	\begin{matrix}  
	0 & 0 \\
	A_{22}V_{021}  & A_{22}V_{022}
	\end{matrix}}
	\end{equation*}
\begin{equation*}
=\paren{
	\begin{matrix}  
	0 & 0 \\
	V_{021}A_{11}  & V_{022}A_{22}
	\end{matrix}}
=V_{0}(	V_{021}A_{11},V_{022}A_{22})=V_{0}(V_{021},V_{022})A,
\end{equation*}
\begin{equation*}
A^2V_{1}(V_{121},V_{122})=
\paren{
	\begin{matrix}  
	A_{11}^2  & 0 \\
	0 & A_{22}^2
	\end{matrix}}
\paren{
	\begin{matrix}  
	0 & 0 \\
	V_{021}  & V_{022}
	\end{matrix}}
=\paren{
	\begin{matrix}  
	0 & 0 \\
	A_{22}^2 V_{021}  & A_{22}^2 V_{022}
	\end{matrix}}
\end{equation*}
\begin{equation*}
=\paren{
	\begin{matrix}  
	0 & 0 \\
	V_{021}A_{11}^2  & V_{022}A_{22}^2 
	\end{matrix}}
=V_{0}(	V_{021}A_{11}^2 ,V_{022}A_{22}^2)=V_{1}(V_{121},V_{122})A^2,
\end{equation*}
\begin{equation*}
A^3V_{2}(V_{0},V_{1},V_{212})=
\paren{
	\begin{matrix}  
	A_{11}^3  & 0 \\
	0 & A_{22}^3
	\end{matrix}}
\paren{
	\begin{matrix}  
	0 & V_{212}\\
	\frac{V_{122}V_{021}}{6}  & \frac{V_{122}V_{022}}{2}
	\end{matrix}}
=\paren{
	\begin{matrix}  
	0 & A_{11}^3V_{212}\\
	\frac{A_{22}^3V_{122}V_{021}}{6}  & \frac{A_{22}^3V_{122}V_{022}}{2}
		\end{matrix}}
\end{equation*}
\begin{equation*}
=\paren{
	\begin{matrix}  
	0 & A_{11}^3V_{212}\\
	\frac{(A_{22}^2V_{122})(A_{22}V_{021})}{6}  & \frac{(A_{22}^2 V_{122})(A_{22}V_{022})}{2}
		\end{matrix}}
=V_{2}(AV_{0},A^2 V_{1}, A_{11}^3 V_{212})=V_{2}(V_{0},V_{1},V_{212})A^3
\mbox{ .}
\end{equation*}

On the other hand, using the Lemma \ref{conmu} we obtain
\begin{equation*}
V_{3}(AV_{0},A^2V_{1},A^3 V_{2})=\sum_{j=1}^{2}jT_{3}^{-1}((A^{j+1}V_{j})(A^{3-j}V_{2-j}))=T_{3}^{-1}((A^2V_{1})^2+2(A^3V_{2})(AV_{0}))
\end{equation*}
\begin{equation*}
=T_{3}^{-1}(A^4(V_{1}^2+2V_{2}V_{0}))=T_{3}^{-1}L_{A}^4(V_{1}^2+2V_{2}V_{0})
=L_{A}^4T_{3}^{-1}(V_{1}^2+2V_{2}V_{0})=A^{4}V_{3}(V_{0},V_{1},V_{2}).
\end{equation*}
Similarly, 
\begin{equation*}
V_{3}(AV_{0},A^2V_{1},A^3 V_{2})=\sum_{j=1}^{2}jT_{3}^{-1}((A^{j+1}V_{j})(A^{3-j}V_{2-j}))=T_{3}^{-1}((A^2V_{1})^2+2(A^3V_{2})(AV_{0}))
\end{equation*}
\begin{equation*}
=T_{3}^{-1}((V_{1}^2+2V_{2}V_{0})A^{4})=T_{3}^{-1}R_{A}^4(V_{1}^2+2V_{2}V_{0})
=R_{A}^4T_{3}^{-1}(V_{1}^2+2V_{2}V_{0})=V_{3}(V_{0},V_{1},V_{2})A^{4}.
\end{equation*}
Assume the claim is true for $3\leq j\leq k-1$ and use again  the Lemma \ref{conmu} to obtain
\begin{equation*}
V_{k}(AV_{0},A^2V_{1},A^3 V_{2})=\sum_{j=1}^{k-1}jT_{k}^{-1}((A^{j+1}V_{j})(A^{k-j}V_{k-1-j}))
\end{equation*}
\begin{equation*}
=\sum_{j=1}^{k-1}jT_{k}^{-1}(A^{k+1}V_{j}V_{k-1-j})=\sum_{j=1}^{k-1}jT_{k}^{-1}L_{A}^{k+1}(V_{j}V_{k-1-j})
=L_{A}^{k+1}\sum_{j=1}^{k-1}jT_{k}^{-1}(V_{j}V_{k-1-j})=A^{k+1}V_{k}(V_{0},V_{1},V_{2}).
\end{equation*}
Similarly,
\begin{equation*}
V_{k}(AV_{0},A^2V_{1},A^3 V_{2})=\sum_{j=1}^{k-1}jT_{k}^{-1}((A^{j+1}V_{j})(A^{k-j}V_{k-1-j}))
\end{equation*}
\begin{equation*}
=\sum_{j=1}^{k-1}jT_{k}^{-1}(V_{j}V_{k-1-j}A^{k+1})=\sum_{j=1}^{k-1}jT_{k}^{-1}R_{A}^{k+1}(V_{j}V_{k-1-j})
=R_{A}^{k+1}\sum_{j=1}^{k-1}jT_{k}^{-1}(V_{j}V_{k-1-j})=V_{k}(V_{0},V_{1},V_{2})A^{k+1}.
\end{equation*}
Thus, the claim follows by induction. \QED
\end{demos}

\begin{coro}\label{v-quasi}
	$V_{k}(\lambda V_{0},\lambda^2 V_{1}, \lambda^3 V_{2})=\lambda^{k+1} V_{k}(V_{0},V_{1},V_{2})$ for every $\lambda\in \complex$, $k\in \natu$, i.e the function  $V_{k}(V_{0},V_{1},V_{2})$ is quasihomogeneous of type $(1,2,3)$ and degree $k+1$.
\end{coro}
\begin{demos}
It is enough to consider $A=\lambda I_{N}$. \QED
\end{demos}

Two remarkable cases of the sequence of functions$\llave{V_{k}(V_{0},V_{1},V_{2})}_{k\in \natu}$ are 
\begin{itemize}
	\item If $V_{-1}=-2I_{N}$, then the equations $V_{-1}V_{0}=V_{-1}V_{1}=0$ implies that $V_{0}=V_{1}=0$. On the other hand, $T_{k}^{-1}(a)=\frac{1}{(k-1)(k+2)}a$ and \eqref{poly} imply that 
	$V_{k}=0$ for $k\geq 2$. In this case we have $V(x)=-\frac{2I_{N}}{x}$.
	 
	\item If $V_{-1}=0$, then $V_{0}$ and $V_{1}$ are arbitrary. Furthermore,  $T_{k}^{-1}(a)=\frac{1}{k(k-1)}a$ and \eqref{poly} imply that 
	
	\begin{equation*}
	V_{k}= \frac{1}{k(k-1)}\sum_{j=1}^{k-1} j  V_{j}V_{k-1-j},
	\end{equation*}
	for $k\geq 2$. In particular, $V_{k}=V_{k}(V_{0},V_{1})$ is a noncommutative  polynomial in the variables $V_{0}$ and $V_{1}$.
\end{itemize}

In the last  case we have an interesting result 
\begin{prop}\label{polyprop}
	\begin{itemize}
		If $V_{-1}=0$, then
		
		\item 	$V_{1}$ is a left divisor of $V_{k}$ for $k\geq 1$. In particular $V_{k}(V_{0},0)=0$ for $k\geq 1$.
		\item  $V_{k}(0,V_{1})=0$ for $k$ even, $V_{k}(0,V_{1})=r_{4k-1}V_{1}^{2k}$ for $k\geq 1$ and $V_{k}(0,V_{1})=r_{4k+1}V_{1}^{2k+1}$ for $k\geq 0$ for some coefficients 
		$r_{4k-1},r_{4k+1}\in [0,1]$. 
	\end{itemize}
	
\end{prop}

\subsection{Polynomial Solutions of the matrix equation $V''(x)=V'(x)V(x)$}

If we want a polynomial solution of degree $\leq n$ for  the equation $V''=V'V$ we have to solve the system of matrix equations 

\begin{equation*}
V_{s}=\frac{1}{s(s-1)}\sum_{j=1}^{s-1} j  V_{j}V_{s-1-j}, \sum_{j=\max\llave{k-1-n,1}}^{n} j  V_{j}V_{k-1-j}=0
\end{equation*}
for $2\leq s \leq n$, $n+1\leq k\leq 2n+1$.

Using the Proposition \ref{polyprop} we have one class of solutions to this problem.

\begin{teor}
	\begin{itemize}
		Let $V_{1}$ to be a nilpotent matrix of degree $n+1\leq N$,
		\item 	If $n=2k$, then $V(x)=\sum_{j=1}^{k}r_{4j-1}V_{1}^{2j}x^j$ is a solution of  $V''=V'V$ for some $\llave{r_{4j-1}}_{1\leq j\leq k}\subset \complex$.
		\item  If $n=2k+1$, then $V(x)=\sum_{j=1}^{k}r_{4j+1}V_{1}^{2j+1}x^j$ is a solution of  $V''=V'V$ for some $\llave{r_{4j+1}}_{1\leq j\leq k}\subset \complex$.
	\end{itemize}
	
\end{teor}

\begin{remark}
In the scalar case we have the integral domain $\complex[x]$, if $V$ is a polynomial such that $\deg(V)\geq 2$ we have that $V''$ is a nonzero polynomial. Applying the function $\deg$ to the equation
	\begin{equation*}
\deg(V'')=\deg(V)-2=\deg(V'V)=\deg(V')+\deg(V)=2\deg(V)-1
	\end{equation*}
	and therefore $\deg(V)=-1$, contradiction. Therefore, $\deg(V)\leq 1$, in the case $\deg(V)=1$ there is no solution of the equation, in fact $V''=0$ and $V'V$ is a nonzero polynomial.
	Thus, we have the trivial constant solution $V(x)=V_{0}$.
\end{remark}

\section{Bispectrality of the Matrix Schr\"{o}dinger Bispectral Operators for polynomial potentials}\label{bisch}

We begin with the definition of the family $\mathcal{P}=\llave{P_{k}}_{k\in \natu}$ which will be used to describe  the map $\theta \mapsto \mathcal{B}$ such that 
$(\psi\mathcal{B})(x,z)=\theta(x)\psi(x,z)$ and the bispectral algebra 
\begin{equation*}
\mathbb{A}=\llave{\theta\in M_{N}(\complex[x])\mid \exists B=B(z,\partial_{z}), (\psi B)(x,z)=\theta(x)\psi(x,z)}.
\end{equation*}
\begin{defin}
	For $k\in \natu$ and $\theta \in M_{N}(\complex[x])$, we define  
	\begin{equation*}
	P_{k}(\theta)=\frac{\theta^{(k)}(0)}{(k-1)!}-\frac{1}{2}\sum_{j=0}^{k}\corch{\frac{\theta^{(j)}(0)}{j!},V_{k-1-j}}.
	\end{equation*}
\end{defin}

Now we study some properties of the sequence $\llave{P_k}_ {k\in \natu}$.

\begin{lemma}[Product Formula for $P_{k}$]\label{product formula}
If $\theta_1, \theta_2 \in M_{N}(\complex[x])$, then
\begin{equation*}\label{product}
P_k (\theta_1 \theta_2)=\sum_{s=0}^{k}\llave{P_{k-s}(\theta_1)\frac{\theta_{2}^{(s)}(0)}{s!}+
\frac{\theta_{1}^{(s)}(0)}{s!}P_{k-s}(\theta_2)}
\end{equation*} 
\end{lemma}
\begin{demos}
By definition, 
\begin{equation*}
P_{k}(\theta_{1} \theta_{2})=
\frac{(\theta_{1}\theta_{2})^{(k)}(0)}{(k-1)!}-\frac{1}{2}\sum_{j=0}^{k}\corch{\frac{(\theta_{1}\theta_{2})^{(j)}(0)}{j!}, V_{k-1-j}}
\end{equation*}	
\begin{equation*}
=\frac{1}{(k-1)!}\sum_{j=0}^{k}\binom{k}{j}\theta_{1}^{(j)}(0)\theta_{2}^{(k-j)}(0)
-\frac{1}{2}\sum_{j=0}^{k}\corch{\frac{1}{j!}\sum_{r=0}^{j}\binom{j}{r}\theta_{1}^{(r)}(0)\theta_{2}^{(j-r)}(0),V_{k-1-j}}
\end{equation*}
\begin{equation*}
=k\sum_{j=0}^{k}\frac{\theta_{1}^{(j)}(0)}{j!}\frac{\theta_{2}^{(k-j)}(0)}{(k-j)!}
-\frac{1}{2}\sum_{j=0}^{k}\sum_{r=0}^{j}\corch{\frac{\theta_{1}^{(r)}(0)}{r!}\frac{\theta_{2}^{(j-r)}(0)}{(j-r)!},V_{k-1-j}}
\end{equation*}
\begin{equation*}
=k\sum_{j=0}^{k}\frac{\theta_{1}^{(j)}(0)}{j!}\frac{\theta_{2}^{(k-j)}(0)}{(k-j)!}
-\frac{1}{2}\sum_{j=0}^{k}\sum_{r=0}^{j}\paren{\corch{\frac{\theta_{1}^{(r)}(0)}{r!}, V_{k-1-j}}\frac{\theta_{2}^{(j-r)}(0)}{(j-r)!}+
\frac{\theta_{1}^{(r)}(0)}{r!} \corch{\frac{\theta_{2}^{(j-r)}(0)}{(j-r)!},V_{k-1-j}}}.
\end{equation*}
However,
\begin{equation*}
P_{k}(\theta)-\frac{\theta^{(k)}(0)}{(k-1)!}=
-\frac{1}{2}\sum_{j=0}^{k}\corch{\frac{\theta^{(j)}(0)}{j!},V_{k-1-j}};
\end{equation*}
for every $\theta\in M_{N}(\complex[x])$.

Therefore, 
\begin{equation*}
\sum_{j=0}^{k}\sum_{r=0}^{j}\paren{\corch{\frac{\theta_{1}^{(r)}(0)}{r!}, V_{k-1-j}}\frac{\theta_{2}^{(j-r)}(0)}{(j-r)!}+
	\frac{\theta_{1}^{(r)}(0)}{r!} \corch{\frac{\theta_{2}^{(j-r)}(0)}{(j-r)!},V_{k-1-j}}}
\end{equation*}
\begin{equation*}
=\sum_{r=0}^{k}\sum_{j=r}^{k}\paren{\corch{\frac{\theta_{1}^{(r)}(0)}{r!}, V_{k-1-j}}\frac{\theta_{2}^{(j-r)}(0)}{(j-r)!}+
	\frac{\theta_{1}^{(j-r)}(0)}{(j-r)!} \corch{\frac{\theta_{2}^{(r)}(0)}{r!},V_{k-1-j}}}
\end{equation*}
\begin{equation*}
=\sum_{r=0}^{k}\sum_{s=0}^{k-r}\paren{\corch{\frac{\theta_{1}^{(r)}(0)}{r!}, V_{k-1-s-r}}\frac{\theta_{2}^{(s)}(0)}{s!}+
	\frac{\theta_{1}^{(s)}(0)}{s!} \corch{\frac{\theta_{2}^{(r)}(0)}{r!},V_{k-1-j}}}
\end{equation*}
\begin{equation*}
=\sum_{s=0}^{k}\sum_{r=0}^{k-s}\paren{\corch{\frac{\theta_{1}^{(r)}(0)}{r!}, V_{k-1-s-r}}\frac{\theta_{2}^{(s)}(0)}{s!}+
	\frac{\theta_{1}^{(s)}(0)}{s!} \corch{\frac{\theta_{2}^{(r)}(0)}{r!},V_{k-1-j}}}
\end{equation*}
\begin{equation*}
=\sum_{s=0}^{k}\llave{\paren{\sum_{r=0}^{k-s}\corch{\frac{\theta_{1}^{(r)}(0)}{r!}, V_{k-1-s-r}}  }\frac{\theta_{2}^{(s)}(0)}{s!}+
	\frac{\theta_{1}^{(s)}(0)}{s!}\paren{\sum_{r=0}^{k-s}\corch{\frac{\theta_{2}^{(r)}(0)}{r!},V_{k-1-j}}} }.
\end{equation*}
This implies that 
\begin{equation*}
P_{k}(\theta_{1}\theta_{2})
=k\sum_{j=0}^{k}\frac{\theta_{1}^{(j)}(0)}{j!}\frac{\theta_{2}^{(k-j)}(0)}{(k-j)!}
\end{equation*}
\begin{equation*}
-\frac{1}{2}\sum_{s=0}^{k}\llave{\paren{\sum_{r=0}^{k-s}\corch{\frac{\theta_{1}^{(r)}(0)}{r!}, V_{k-1-s-r}}  }\frac{\theta_{2}^{(s)}(0)}{s!}+
	\frac{\theta_{1}^{(s)}(0)}{s!}\paren{\sum_{r=0}^{k-s}\corch{\frac{\theta_{2}^{(r)}(0)}{r!},V_{k-1-j}}} }
\end{equation*}
\begin{equation*}
=k\sum_{j=0}^{k}\frac{\theta_{1}^{(j)}(0)}{j!}\frac{\theta_{2}^{(k-j)}(0)}{(k-j)!}
\end{equation*}
\begin{equation*}
+\sum_{s=0}^{k}\llave{\paren{P_{k-s}(\theta_{1})-\frac{\theta_{1}^{(k-s)}(0)}{(k-s-1)!}}\frac{\theta_{2}^{(s)}(0)}{s!}+
	\frac{\theta_{1}^{(s)}(0)}{s!}\paren{P_{k-s}(\theta_{2})-\frac{\theta_{2}^{(k-s)}(0)}{(k-s-1)!}} }
\end{equation*}
\begin{equation*}
=k\sum_{j=0}^{k}\frac{\theta_{1}^{(j)}(0)}{j!}\frac{\theta_{2}^{(k-j)}(0)}{(k-j)!}
\end{equation*}
\begin{equation*}
+\sum_{s=0}^{k}P_{k-s}(\theta_{1})\frac{\theta_{2}^{(s)}(0)}{s!}-\sum_{s=0}^{k}\frac{\theta_{1}^{(k-s)}(0)}{(k-s-1)!}\frac{\theta_{2}^{(s)}(0)}{s!}
\end{equation*}
\begin{equation*}
+\sum_{s=0}^{k}	\frac{\theta_{1}^{(s)}(0)}{s!}P_{k-s}(\theta_{2})-\sum_{s=0}^{k}\frac{\theta_{1}^{(s)}(0)}{s!}\frac{\theta_{1}^{(k-s)}(0)}{(k-s-1)!} 
\end{equation*}
\begin{equation*}
=k\sum_{j=0}^{k}\frac{\theta_{1}^{(j)}(0)}{j!}\frac{\theta_{2}^{(k-j)}(0)}{(k-j)!}
+\sum_{s=0}^{k}\llave{P_{k-s}(\theta_{1})\frac{\theta_{2}^{(s)}(0)}{s!}-\frac{\theta_{1}^{(s)}(0)}{s!}P_{k-s}(\theta_{2})}
\end{equation*}
\begin{equation*}
-\sum_{s=0}^{k}(k-s)\frac{\theta_{1}^{(k-s)}(0)}{(k-s)!}\frac{\theta_{2}^{(s)}(0)}{s!}
-\sum_{s=0}^{k}(k-s)\frac{\theta_{1}^{(s)}(0)}{s!}\frac{\theta_{1}^{(k-s)}(0)}{(k-s)!} 
\end{equation*}
\begin{equation*}
=k\sum_{j=0}^{k}\frac{\theta_{1}^{(j)}(0)}{j!}\frac{\theta_{2}^{(k-j)}(0)}{(k-j)!}
+\sum_{s=0}^{k}\llave{P_{k-s}(\theta_{1})\frac{\theta_{2}^{(s)}(0)}{s!}-\frac{\theta_{1}^{(s)}(0)}{s!}P_{k-s}(\theta_{2})}
\end{equation*}
\begin{equation*}
-\sum_{s=0}^{k}s\frac{\theta_{1}^{(s)}(0)}{s!}\frac{\theta_{2}^{(k-s)}(0)}{(k-s)!}
-\sum_{s=0}^{k}(k-s)\frac{\theta_{1}^{(s)}(0)}{s!}\frac{\theta_{1}^{(k-s)}(0)}{(k-s)!} 
\end{equation*}
\begin{equation*}
=\sum_{s=0}^{k}\llave{P_{k-s}(\theta_{1})\frac{\theta_{2}^{(s)}(0)}{s!}-\frac{\theta_{1}^{(s)}(0)}{s!}P_{k-s}(\theta_{2})}.
\end{equation*}
Thus, 
\begin{equation*}
P_k (\theta_1 \theta_2)=\sum_{s=0}^{k}\llave{P_{k-s}(\theta_1)\frac{\theta_{2}^{(s)}(0)}{s!}+
	\frac{\theta_{1}^{(s)}(0)}{s!}P_{k-s}(\theta_2)}.
\end{equation*} \QED
\end{demos}

\begin{remark}
	If $V_{j}=0$ for every $j\in \natu$ the product formula specializes into the Leibniz rule
	\begin{equation*}
	(\theta_{1}\theta_{2})^{(k)}(x)=\sum_{s=0}^{k}\binom{k}{s}\theta_{1}^{(k-s)}(x)\theta_{2}^{(s)}(x).
	\end{equation*}
In this case, $P_{k}(\theta)=\frac{\theta^{(k)}(0)}{(k-1)!}=k\frac{\theta^{(k)}(0)}{k!}$ and applying the Product Formula \eqref{product formula} turns out 
\begin{equation*}
P_{k}(\theta_{1}\theta_{2})=\sum_{s=0}^{k}\llave{(k-s)\frac{\theta_{1}^{(k-s)}(0)}{(k-s)!}\frac{\theta_{2}^{(s)}(0)}{s!}+
\frac{\theta_{1}^{(s)}(0)}{s!}(k-s)\frac{\theta_{2}^{(k-s)}(0)}{(k-s)!}}
\end{equation*}
\begin{equation*}
=\sum_{s=0}^{k}\llave{(k-s)\frac{\theta_{1}^{(k-s)}(0)}{(k-s)!}\frac{\theta_{2}^{(s)}(0)}{s!}+
	s\frac{\theta_{1}^{(k-s)}(0)}{(k-s)!}\frac{\theta_{2}^{(s)}(0)}{s!}}
\end{equation*}
\begin{equation*}
=k\sum_{s=0}^{k}\frac{\theta_{1}^{(k-s)}(0)}{(k-s)!}\frac{\theta_{2}^{(s)}(0)}{s!}
=k\frac{(\theta_{1}\theta_{2})^{(k)}(0)}{k!},
\end{equation*}
in the words $(\theta_{1}\theta_{2})^{(k)}(0)=\sum_{s=0}^{k}\binom{k}{s}\theta_{1}^{(k-s)}(0)\theta_{2}^{(s)}(0)$. Since $\theta_{1}$ and $\theta_{2}$ were arbitrary we can change them by their translations $\theta_{1}(x+\cdot)$ and $\theta_{2}(x+\cdot)$ to obtain 	
$(\theta_{1}\theta_{2})^{(k)}(x)=\sum_{s=0}^{k}\binom{k}{s}\theta_{1}^{(k-s)}(x)\theta_{2}^{(s)}(x)$, i.e., the Leibniz rule.
	\end{remark}

If we consider the formal power series $V(x)=\sum_{j=-1}^{\infty}V_{j}x^j$ we can write the family in a nice form as stated in the following theorem.

\begin{teor}
	For every $k\in \natu$ we have 
	\begin{equation*}
	P_{k}=\frac{1}{k!}\frac{d^k}{dx^k} \Big|_{x=0}\paren{kI+\frac{1}{2}x\; ad(V)}
	\end{equation*}
\end{teor}
\begin{demos}
	Since $V(x)=\sum_{j=-1}^{\infty}V_{j}x^j$ we have $xV(x)=\sum_{j=-1}^{\infty}V_{j}x^{j+1}=\sum_{l=0}^{\infty}V_{l-1}x^l$ and $V_{l-1}=
	\frac{1}{l!} 	\frac{d^l}{dx^l} \Big|_{x=0} (xV(x))$. Therefore, 
\begin{equation*}
\sum_{j=0}^{k}\corch{\frac{\theta^{(j)}(0)}{j!},V_{k-1-j}}=
\sum_{j=0}^{k}\corch{\frac{\theta^{(j)}(0)}{j!},\frac{1}{(k-j)!} \frac{d^{k-j}}{dx^{k-j}}\Big|_{x=0}(xV(x))}
=\frac{1}{k!}\frac{d^k}{dx^k}\Big|_{x=0}  (x\corch{\theta, V}).
\end{equation*}
Thus, 
\begin{equation*}
P_{k}(\theta)=\frac{\theta^{(k)}(0)}{(k-1)!}-\frac{1}{2}\sum_{j=0}^{k}\corch{\frac{\theta^{(j)}(0)}{j!},V_{k-1-j}}
=\frac{\theta^{(k)}(0)}{(k-1)!}-\frac{1}{2}\frac{1}{k!}\frac{d^k}{dx^k}\Big|_{x=0}(x\corch{\theta, V})
\end{equation*}
\begin{equation*}
=\frac{\theta^{(k)}(0)}{(k-1)!}+\frac{1}{2}\frac{1}{k!}\frac{d^k}{dx^k}\Big|_{x=0}(x ad(V)(\theta))
=\frac{1}{k!}\frac{d^k}{dx^k}\Big|_{x=0} \paren{kI+\frac{1}{2}xad(V)}(\theta).
\end{equation*}
Since $\theta$ is arbitrary we have the assertion. \QED
\end{demos}

\begin{coro}
For every $k\in \natu$,	$P_{k}(V)=\frac{V^{(k-1)}(0)}{(k-1)!}$ .
\end{coro}

\begin{defin}
	For $m\in \natu$ define 
		\begin{equation*}
	A_{1}^{[m]}=
	\paren{
		\begin{matrix} 
		\frac{V_0}{2} & \frac{1}{2}V_{-1}+I_N & 0 &   \cdots & 0 & 0 & 0  \\
		\frac{V_1}{2} & \frac{V_0}{2} & \frac{1}{2}V_{-1}+2I_{N} &  \cdots & 0 & 0 & 0 \\
		. & . & \cdots & .  & \cdots & . & . \\
		. & . & \cdots & . & \cdots & . & . \\
		. & . & \cdots & . & \cdots & . & . \\
		\frac{V_{m-2}}{2}  &  	\frac{V_{m-3}}{2}  & 	\frac{V_{m-4}}{2}  &  \cdots & \frac{V_0}{2} & \frac{1}{2}V_{-1}+(m-1)I_{N}  & 0 \\
		\frac{V_{m-1}}{2}  & \frac{V_{m-2}}{2}  & \frac{V_{m-3}}{2}   &  \cdots &  \frac{V_1}{2}  & \frac{V_0}{2} & \frac{1}{2}V_{-1}+mI_{N} \\
		\frac{V_m}{2} & \frac{V_{m-1}}{2} & \frac{V_{m-2}}{2}   &  \cdots & \frac{V_2}{2} & \frac{V_1}{2} & \frac{V_0}{2}\\
		\end{matrix}},
	\end{equation*}
	\begin{equation*}
	A_{2}^{[m]}=	\paren{
		\begin{matrix} 
		V_{m+1} & V_{m}& \cdots & V_{1} \\
		V_{m+2} & V_{m+1}  & \cdots & V_{2}\\
		. & . & \cdots & . \\
		. & . & \cdots & .\\
		. & . & \cdots & .\\
		\end{matrix}},
	\end{equation*}
	and for $\theta \in M_{N}(\complex[x])$ we define $P_{1}^{m+1}(\theta)=(P_{1}(\theta), P_{2}(\theta), \cdots , P_{m}(\theta), P_{m+1}(\theta))^{T}$ and 
\\	$P_{m+2}^{\infty}(\theta)=(P_{m+2}(\theta), P_{m+3}(\theta), \cdots )^{T}$. 
\end{defin}

Note that $A_{1}^{[m]}$, $A_{2}^{[m]}$ depend on $V$ and $P_{1}^{m+1}(\theta)$, 	$P_{m+2}^{\infty}(\theta)$ depend on $\theta$. The following lemma  gives a simplification of these  matrices for $m$ large enough when $V$ is a polynomial.

\begin{lemma}
	If $V$ is a polynomial of degree $n$ and $m+1=nq+r$, $q\geq 1$, $ 0\leq r <n$, then 
{	\scriptsize
	\begin{equation*}
A_{1}^{[m]}=
\paren{
	\begin{matrix} 
	A_{1}^{[r-1]} & r T_{r}^{r-1} \hspace{0.1cm},  0_{r\times (n-r)} & 0_{n\times n} &   \cdots & 0_{n\times n}  & 0_{n\times n}  & 0_{n\times n}   \\
	\frac{1}{2}(A_{2}^{[n-1]})^{1,2,\cdots,r} & A_{1}^{[n-1]}+rS_{n} & (n+r)T_{n}^{n-1} &  \cdots & 0_{n\times n} & 0_{n\times n}  & 0_{n\times n}  \\
	 0_{n\times n}  & 	\frac{1}{2}A_{2}^{[n-1]} & A_{1}^{[n-1]}+(n+r)S_{n} &  \cdots & 0_{n\times n} & 0_{n\times n}  & 0_{n\times n}  \\
	. & . & \cdots & .  & \cdots & . & . \\
	. & . & \cdots & . & \cdots & . & . \\
	. & . & \cdots & . & \cdots & . & . \\
0_{n\times n}  &0_{n\times n}  & 0_{n\times n}   &  \cdots &  \frac{1}{2}A_{2}^{[n-1]}  &A_{1}^{[n-1]}+(n(q-2)+r)S_{n}& (n(q-1)+r)T_{n}^{n-1}\\
0_{n\times n}   & 0_{n\times n}  & 0_{n\times n}   &  \cdots &0_{n\times n}   & \frac{1}{2}A_{2}^{[n-1]} & A_{1}^{[n-1]}+(n(q-1)+r)S_{n}\\
	\end{matrix}},
\end{equation*}
}
and 
{	\footnotesize
	\begin{equation*}
	A_{2}^{[m]}=
	\paren{
		\begin{matrix} 
0_{n \times (m+1-n)} 	 & A_{2}^{[n-1]} \\
0_{\infty \times (m+1-n)} 	 & 0_{\infty \times n} 	\\
		\end{matrix}}\mbox{ ,}
	\end{equation*} 
}
i.e., $A_{1}^{[m]}$ is a block tridiagonal matrix and $	A_{2}^{[m]}$ is a block upper triangular matrix.
\end{lemma}
\begin{demos}
	If $V$ is a polynomial of degree $n$ the assertion about $	A_{2}^{[m]}$ is clear. On the other hand, note that we can write 
	\begin{equation}\label{partition}
		A_{1}^{[m]}=
	\paren{
		\begin{matrix} 
		A_{1}^{[m-n]}	 & 	\begin{matrix} 
		0_{(m+1-2n) \times n } 	\\ \\
	(m+1-n)T_{n}^{n-1} 	 	\\
		\end{matrix}\\ \\
		
\begin{matrix} 
0_{n \times (m+1-2n) } 	& 
\frac{1}{2}A_{2}^{[n-1]}	\\
\end{matrix}	 &  	A_{1}^{[n-1]}+(m+1-n)S_{n}	\\
		\end{matrix}}
	\end{equation}
	Since the 	$A_{1}^{[m]}$ is a block matrix of size $(q+1) \times (q+1)$ we can use induction over $q$. Notice that the assertion is clear for $q=1$ because in this case $m+1=n+r$ and 
	\begin{equation*}
	A_{1}^{[m]}=
	\paren{
		\begin{matrix} 
		A_{1}^{[r-1]} & r T_{r}^{r-1} \hspace{0.1cm},  0_{r\times (n-r)}   \\
		\frac{1}{2}(A_{2}^{[n-1]})^{1,2,\cdots,r} & A_{1}^{[n-1]}+rS_{n}  \\
		\end{matrix}}.
	\end{equation*}
Now let $m+1=nq+r$ and assume the assertion for $m-n$ or equivalently for a matrix of size $q\times q$. Therefore, 
	 
	 {	\scriptsize
	 	\begin{equation*}
	 	A_{1}^{[m-n]}=
	 	\paren{
	 		\begin{matrix} 
	 		A_{1}^{[r-1]} & r T_{r}^{r-1} \hspace{0.1cm},  0_{r\times (n-r)} & 0_{n\times n} &   \cdots & 0_{n\times n}  & 0_{n\times n}  & 0_{n\times n}   \\
	 		\frac{1}{2}(A_{2}^{[n-1]})^{1,2,\cdots,r} & A_{1}^{[n-1]}+rS_{n} & (n+r)T_{n}^{n-1} &  \cdots & 0_{n\times n} & 0_{n\times n}  & 0_{n\times n}  \\
	 		0_{n\times n}  & 	\frac{1}{2}A_{2}^{[n-1]} & A_{1}^{[n-1]}+(n+r)S_{n} &  \cdots & 0_{n\times n} & 0_{n\times n}  & 0_{n\times n}  \\
	 		. & . & \cdots & .  & \cdots & . & . \\
	 		. & . & \cdots & . & \cdots & . & . \\
	 		. & . & \cdots & . & \cdots & . & . \\
	 		0_{n\times n}  &0_{n\times n}  & 0_{n\times n}   &  \cdots &  \frac{1}{2}A_{2}^{[n-1]}  &A_{1}^{[n-1]}+(n(q-3)+r)S_{n}& (n(q-2)+r)T_{n}^{n-1}\\
	 		0_{n\times n}   & 0_{n\times n}  & 0_{n\times n}   &  \cdots &0_{n\times n}   & \frac{1}{2}A_{2}^{[n-1]} & A_{1}^{[n-1]}+(n(q-2)+r)S_{n}\\
	 		\end{matrix}}.
	 	\end{equation*}
	 }
 If we replace this in \eqref{partition} we obtain the claim. Thus, the assertion follows by induction.\QED
\end{demos}

The following theorem characterizes bispectrality using the family $\llave{P_{k}}_{k\in \natu}$.

\begin{teor}\label{bispectralth}
Let {\small $$ \Gamma= \left\lbrace\theta\in M_{N}(\complex[x])\mid P_{0}(\theta)=0, V_{-1}e_{1}(A_{1}^{[m]})^{k}P_{1}^{m+1}(\theta)=0,
	 A_{2}^{[m]}(A_{1}^{[m]})^{k}P_{1}^{m+1}(\theta)=0, k\geq 0, P_{m+2}^{\infty}(\theta)=0, m=\deg(\theta) \right\rbrace$$} then 
$\Gamma=\mathbb{A}$. Moreover, for each $\theta$ we have an explicit expression for the operator $B$ such that 
\begin{equation*}
(\psi B)(x,z)=\theta(x)\psi(x,z).
\end{equation*}
\end{teor}
\begin{remark}
	Before proving the Theorem \ref{bispectralth} we observe that since $A_{1}^{[m]}\in M_{(m+1)N}(\complex)$ the Cayley-Hamilton Theorem implies that we can assume that \\
	{\small $\Gamma= \left\lbrace\theta\in M_{N}(\complex[x])\mid P_{0}(\theta)=0, V_{-1}e_{1}(A_{1}^{[m]})^{k}P_{1}^{m+1}(\theta)=0,
	A_{2}^{[m]}(A_{1}^{[m]})^{k}P_{1}^{m+1}(\theta)=0, 0\leq k\leq (m+1)N-1,\right.$\\  $ \left. P_{m+2}^{\infty}(\theta)=0, m=\deg(\theta) \right\rbrace$.}
\end{remark}
\begin{demos}
	If we consider $\theta(x)=\sum_{j=0}^{m}a_{j}x^j$ and $B(z,\partial_{z})=\sum_{j=0}^{m}\partial_{z}^{j}\cdot b_{j}(z)$ then, 
\begin{equation*}
\Lambda(x,z)=e^{-xz}((\psi B)(x,z)-\theta(x)\psi(x,z))
\end{equation*}
\begin{equation*}
=e^{-xz}\paren{\sum_{j=0}^{m}\partial_{z}^{j}\paren{\paren{Iz+\frac{1}{2}V(x)}e^{xz}}\cdot b_{j}(z)-\sum_{j=0}^{m}a_{j}x^j\paren{Iz+\frac{1}{2}V(x)}e^{xz}}
\end{equation*}
\begin{equation*}
=e^{-xz}\paren{\sum_{j=0}^{m}\paren{\sum_{l=0}^{j}\binom{j}{l}\partial_{z}^l \paren{Iz+\frac{1}{2}V(x)}\partial_{z}^{j-l}(e^{xz})}b_{j}(z) 
-\sum_{j=0}^{m}a_{j}x^j z e^{xz}-\sum_{j=0}^{m}\frac{a_{j}}{2}x^j V(x) e^{xz} }
\end{equation*}
\begin{equation*}
=e^{-xz}\paren{\sum_{j=0}^{m}\paren{\paren{Iz+\frac{1}{2}V(x)}x^j e^{xz}+jx^{j-1}e^{xz}}b_{j}(z) 
	-\sum_{j=0}^{m}a_{j}x^j z e^{xz}-\sum_{j=0}^{m}\frac{a_{j}}{2}x^j V(x) e^{xz} }
\end{equation*}
\begin{equation*}
=\sum_{j=0}^{m}b_{j}(z)x^j z +\sum_{j=0}^{m}\frac{1}{2}V(x)x^j b_{j}(z)+ \sum_{j=0}^{m}j x^{j-1} b_{j}(z)
-\sum_{j=0}^{m}a_{j}x^j z -\sum_{j=0}^{m}\frac{a_{j}}{2}x^j V(x). 
\end{equation*}
Writing $V(x)=\sum_{j=-1}^{\infty}V_{j}x^{j}$ then, 
\begin{equation*}
\Lambda(x,z)=\sum_{j=0}^{m}b_{j}(z)x^j z +\sum_{j=0}^{m}\frac{1}{2}\paren{\sum_{k=-1}^{\infty}V_{j}x^{j}}x^j b_{j}(z)+ \sum_{j=0}^{m}j x^{j-1} b_{j}(z)
-\sum_{j=0}^{m}a_{j}x^j z -\sum_{j=0}^{m}\frac{a_{j}}{2}x^j\paren{\sum_{k=-1}^{\infty}V_{j}x^{j}}
\end{equation*}
\begin{equation*}
=\sum_{j=0}^{m}b_{j}(z)x^j z +\sum_{j=0}^{m}\sum_{k=-1}^{\infty}\frac{1}{2}V_{k} b_{j}(z) x^{k+j}+ \sum_{j=0}^{m}j x^{j-1} b_{j}(z)
-\sum_{j=0}^{m}a_{j}x^j z -\sum_{j=0}^{m}\sum_{k=-1}^{\infty}\frac{a_{j}}{2}V_{k}x^{k+j}.
\end{equation*}
Let $s=k+j$ then $s$ varies from $-1$ to $\infty$. 
\begin{equation*}
\Lambda(x,z)=\sum_{j=0}^{m}b_{j}(z)x^j z +\sum_{s=-1}^{\infty} \sum_{j=0}^{s+1}\frac{1}{2}V_{s-j}b_{j}(z)x^s
+\sum_{j=1}^{m}jx^{j-1} b_{j}(z) - \sum_{j=0}^{m}a_{j}x^j z - \sum_{s=-1}^{\infty} \sum_{j=0}^{s+1}\frac{1}{2}a_{j}V_{s-j}x^s
\end{equation*}
\begin{equation*}
=\sum_{j=0}^{m}b_{j}(z)z x^j +\sum_{s=-1}^{\infty} \paren{\sum_{j=0}^{s+1}\frac{1}{2}V_{s-j}b_{j}(z)}x^s
+\sum_{s=0}^{m-1}(s+1)b_{s+1}(z)x^s - \sum_{j=0}^{m}a_{j}z x^j - \sum_{s=-1}^{\infty} \paren{\sum_{j=0}^{s+1}\frac{1}{2}a_{j}V_{s-j}}x^s
\end{equation*}

\begin{equation*}
=\frac{1}{2}(V_{-1}b_{0}(z)-a_{0}V_{-1}) + 
\sum_{s=0}^{m-1}\paren{\sum_{j=0}^{s+1}\frac{1}{2}V_{s-j}b_{j}(z)+ b_{s}(z)z + (s+1)b_{s+1}(z)- a_{s}z - \sum_{j=0}^{s+1}\frac{1}{2}a_{j}V_{s-j}}x^s
\end{equation*}
\begin{equation*}
+\paren{b_{m}(z)z+\sum_{j=0}^{m+1}\frac{1}{2}V_{m-j}b_{j}(z)-a_{m}z- \sum_{j=0}^{m+1}\frac{1}{2}a_{j}V_{m-j}}x^{m}
\end{equation*}
\begin{equation*}
+\sum_{s=m+1}^{\infty}\paren{\sum_{j=0}^{s+1}\frac{1}{2}(V_{s-j}b_{j}(z)-a_{j}V_{s-j})}x^s
\end{equation*}

	 if, and only if, 
	\begin{equation*}
	V_{-1}b_{0}-a_{0}V_{-1}=0,
	\end{equation*}
	\begin{equation*}
	(b_{s}(z)-a_{s})z+(s+1)b_{s+1}(z)+\frac{1}{2}\sum_{k=0}^{s}(V_{s-k}b_{k}-a_{k}V_{s-k})+\frac{1}{2}(V_{-1}b_{s+1}-a_{s+1}V _{-1})=0,
	\end{equation*}
	for $0\leq s \leq m-1$.
	\begin{equation*}
	(b_{m}(z)-a_{m})z+\frac{1}{2}\sum_{k=0}^{m}(V_{m -k}b_{k}-a_{k}V_{m-k})=0,
	\end{equation*}
	\begin{equation*}
	\sum_{k=0}^{m}(V_{s-k}b_k - a_k V_{s-k})=0,
	\end{equation*}
	for $s\geq m+1$. 

If we define $c_{j}(z)=b_{j}(z)-a_{j}$ we have 
\begin{equation*}
\frac{1}{2}V_{-1}c_{0}(z)=-P_{0}(\theta),
\end{equation*}
{\footnotesize
\begin{equation*}
\paren{
	\begin{matrix} 
	z+\frac{V_0}{2} & \frac{1}{2}V_{-1}+I_N & 0 &   \cdots & 0 & 0 & 0  \\
	\frac{V_1}{2} & z+\frac{V_0}{2} & \frac{1}{2}V_{-1}+2I_{N} &  \cdots & 0 & 0 & 0 \\
	. & . & \cdots & .  & \cdots & . & . \\
	. & . & \cdots & . & \cdots & . & . \\
	. & . & \cdots & . & \cdots & . & . \\
	\frac{V_{m-2}}{2}  &  	\frac{V_{m-3}}{2}  & 	\frac{V_{m-4}}{2}  &  \cdots & z+\frac{V_0}{2} & \frac{1}{2}V_{-1}+(m-1)I_{N}  & 0 \\
	\frac{V_{m-1}}{2}  & \frac{V_{m-2}}{2}  & \frac{V_{m-3}}{2}   &  \cdots &  \frac{V_1}{2}  & z+\frac{V_0}{2} & \frac{1}{2}V_{-1}+mI_{N} \\
	\frac{V_m}{2} & \frac{V_{m-1}}{2} & \frac{V_{m-2}}{2}   &  \cdots & \frac{V_2}{2} & \frac{V_1}{2} & z+\frac{V_0}{2}\\
	\end{matrix}}
\paren{
	\begin{matrix} 
	c_{0}(z)  \\
	c_{1}(z) \\
	. \\
	. \\
	. \\
	c_{m-2}(z)  \\
	c_{m-1}(z)\\
	c_m (z)\\
	\end{matrix}}
=\paren{
	\begin{matrix} 
-P_{1}(\theta)\\
-	P_{2}(\theta)\\
	. \\
	. \\
	. \\
-	P_{m-1}(\theta) \\
-	P_{m}(\theta)\\
 -   P_{m+1}(\theta)\\
	\end{matrix}}
\end{equation*}}

and 

\begin{equation*}
\paren{
	\begin{matrix} 
	V_{m+1} & V_{m}& \cdots & V_{1} \\
	V_{m+2} & V_{m+1}  & \cdots & V_{2}\\
	. & . & \cdots & . \\
	. & . & \cdots & .\\
	. & . & \cdots & .\\
	\end{matrix}}
\paren{
	\begin{matrix} 
	c_{0}(z) \\
	c_{1}(z) \\
	. \\
	. \\
	. \\
	c_{m-2}(z)\\
	c_{m-1}(z)\\
	c_m (z)\\
	\end{matrix}}=
\paren{
	\begin{matrix} 
-	P_{m+2}(\theta)\\
-	P_{m+3}(\theta)\\
	. \\
	. \\
	. \\
	\end{matrix}}.
\end{equation*}
Using the notation defined above 
\begin{equation*}
\frac{1}{2}V_{-1}c_{0}(z)=-P_{0}(\theta), (A_{1}^{[m]}+z)c(z)=-P_{1}^{m+1}(\theta), A_{2}^{[m]}c(z)=-P_{m+2}^{\infty}(\theta). 
\end{equation*}

However, $(A_{1}^{[m]}+z)^{-1}=\sum_{k=0}^{\infty}\frac{(-A_{1}^{[m]})^{k}}{z^{k+1}}$ implies that $c(z)=- \sum_{k=0}^{\infty}\frac{(-A_{1}^{[m]})^{k}}{z^{k+1}}P_{1}^{m+1}(\theta)$ 
and $P_{m+2}^{\infty}(\theta)=A_{2}^{[m]}\sum_{k=0}^{\infty}\frac{(-A_{1}^{[m]})^{k}}{z^{k+1}}P_{1}^{m+1}(\theta)$, using $z$ as variable we obtain $P_{0}(\theta)=0$,
$A_{2}^{[m]}(A_{1}^{[m]})^{k}P_{1}^{m+1}(\theta)=0, k\geq 0$ and $P_{m+2}^{\infty}(\theta)=0$. Furthermore, $c_{s}(z)=-e_{s+1}(A_{1}^{[m]}+z)^{-1}P_{1}^{m+1}(\theta)$ for $0\leq s \leq m$. In particular $c_{0}(z)=
-e_{1}(A_{1}^{[m]}+z)^{-1}P_{1}^{m+1}(\theta)$ then, $V_{-1}e_{1}(A_{1}^{[m]}+z)^{-1}P_{1}^{m+1}(\theta)=\sum_{k=0}^{\infty}V_{-1}e_{1}\frac{(-A_{1}^{[m]})^{k}}{z^{k+1}}P_{1}^{m+1}(\theta)=0$. Using $z$ as  variable we obtain
$V_{-1}e_{1}(A_{1}^{[m]})^{k}P_{1}^{m+1}(\theta)=0$ for every $k\in \natu$.
We shall now use this remark to  
conclude the proof of the theorem. 

If $\theta \in \mathbb{A}$, then there exists $B=\sum_{j=0}^{m}\partial_{z}^j \cdot b_{j}(z) $ such that 
\begin{equation*}
\Lambda(x,z)=e^{-xz}((\psi B)(x,z)-\theta(x)\psi(x,z))=0. 
\end{equation*}
But this is equivalent to $P_{0}(\theta)=0$, $V_{-1}e_{1}(A_{1}^{[m]}+z)^{-1}P_{1}^{m+1}(\theta)=0$,$(A_{1}^{[m]}+z)c(z)=P_{1}^{m+1}(\theta)$, $A_{2}^{[m]}c(z)=P_{m+2}^{\infty}(\theta)$ with $c(z)=b(z)-a$, $b=(b_0, \cdots ,b_m)$, 
$a=(a_{0}, \cdots , a_{m})$. 

By the previous remark we have 
\begin{equation*}
P_{0}(\theta)=0, V_{-1}e_{1}(A_{1}^{[m]}+z)^{-1}P_{1}^{m+1}(\theta)=0, A_{2}^{[m]}(A_{1}^{[m]})^{k}P_{1}^{m+1}(\theta)=0, k\geq 0 \; \text{and}\;  P_{m+2}^{\infty}(\theta)=0. 
\end{equation*}
Then $\theta \in \Gamma$. Since $\theta \in \mathbb{A}$ was arbitrary we have  $\mathbb{A}\subset \Gamma$. 

On the other hand, if $\theta\in \Gamma$, then $P_{0}(\theta)=0$, $V_{-1}e_{1}(A_{1}^{[m]}+z)^{-1}P_{1}^{m+1}(\theta)=0$, $A_{2}^{[m]}(A_{1}^{[m]})^{k}P_{1}^{m+1}(\theta)=0, k\geq 0$ and $P_{m+2}^{\infty}(\theta)=0$. 

Taking 
\begin{equation*}
b_{j}(z)=a_{j}+e_{j}(A_{1}^{[m]}+z)^{-1}P_{1}^{m+1}(\theta),
\end{equation*}
for $0\leq j \leq m$.

We have  $c(z)= -\sum_{k=0}^{\infty}\frac{(-A_{1}^{[m]})^{k}}{z^{k+1}}P_{1}^{m+1}(\theta)=(A_{1}^{[m]}+z)^{-1}P_{1}^{m+1}(\theta)$ and therefore
\begin{equation*}
\frac{1}{2}V_{-1}c_{0}(z)=-P_{0}(\theta),  (A_{1}^{[m]}+z)c(z)=-P_{1}^{m+1}(\theta), A_{2}^{[m]}c(z)=-P_{m+2}^{\infty}(\theta).
\end{equation*}
 By the previous arguments we obtain that 
\begin{equation*}
\Lambda(x,z)=e^{-xz}((\psi B)(x,z)-\theta(x)\psi(x,z))=0,
\end{equation*}
with $B=\sum_{j=0}^{m}b_{j}(z)\cdot \partial_{z}^j$. This implies that $\theta \in \mathbb{A}$. Since $\theta \in \Gamma$ was arbitrary we have  $\Gamma \subset \mathbb{A}$. 

Thus, $\Gamma=\mathbb{A}$ and for every $\theta \in \mathbb{A}$ there exists a unique operator $B=\sum_{j=0}^{m}\partial_{z}^j \cdot b_{j}(z) $ given by 
\begin{equation*}
b_{j}(z)=a_{j}-e_{j}(A_{1}^{[m]}+z)^{-1}P_{1}^{m+1}(\theta),
\end{equation*}
for $0\leq j \leq m$, such that 
\begin{equation*}
(\psi B)(x,z)=\theta(x)\psi(x,z). 
\end{equation*}
This concludes the proof of the assertion.
\end{demos}

\begin{coro}
	For $\theta \in \mathbb{A}$ the operator  $B=\sum_{j=0}^{m}\partial_{z}^j \cdot b_{j}(z) $  such that
	\begin{equation*}
	(\psi B)(x,z)=\theta(x)\psi(x,z).
	\end{equation*}
   satisfies $\lim\limits_{z \rightarrow \infty}b_j (z)= a_{j}$ for $0\leq j \leq m$.
\end{coro}

In the following result we rewrite the expressions defining the algebra $\Gamma$ for another more simple to remind.

\begin{lemma}\label{description}
	The algebra $\Gamma$ is exactly the set of all polynomial $\theta\in M_{N}(\complex)[x]$, $m=\deg(\theta)$  such that 
	$\corch{\theta,V}$ is a polynomial of degree $\leq m$ and 
	\begin{equation*}
	V_{-1}e_{1}(A_{1}^{[m]})^{k} P_{1}^{m+1}(\theta)=0, 
 A_{2}^{[m]}(A_{1}^{[m]})^{k}P_{1}^{m+1}(\theta)=0,
	\end{equation*}
	 
	for $0\leq k \leq (m+1)N-1$.
\end{lemma}

\begin{demos}
Note that  $k\geq m+1$ implies  $P_{k}(\theta)=\frac{\theta^{(k)}(0)}{(k-1)!}+\frac{1}{2}\frac{1}{k!}\frac{d^k}{dx^k}\Big|_{x=0}(x ad(V)(\theta))
=\frac{1}{2}\frac{1}{k!}\frac{d^k}{dx^k}\Big|_{x=0}(x ad(V)(\theta)).$ 

Since $\frac{d^k}{dx^k}\Big|_{x=0}(x ad(V)(\theta))=P_{k}(\theta)=0$ for $k\geq m+2$ we have that $x\corch{\theta,V}$ is a polynomial of degree $\leq m+1$.
Furthermore, $P_{0}(\theta)=\frac{1}{2}(x ad(V))(\theta)\Big|_{x=0}=0$ we have that $x\corch{\theta,V}\Big|_{x=0}=0$. However, since $x\corch{\theta,V}$ is a polynomial we have that 
$\corch{\theta,V}$ is a polynomial of degree $\leq m$. Moreover, we have the restrictions $	V_{-1}e_{1}(A_{1}^{[m]})^{k} P_{1}^{m+1}(\theta)=0$, $A_{2}^{[m]}(A_{1}^{[m]})^{k}P_{1}^{m+1}(\theta)=0$, for $0\leq k \leq (m+1)N-1$. \QED
\end{demos}

Now we try to find some solutions of these equations. To do this we put restrictions on the matrix $V_0$ and $V_{1}$ to obtain $V\in \mathbb{A}$. 
We begin with a definition

\begin{defin}
	We define the grading $\deg_{1,2}$ on the ring $\complex \langle V_0 ,V_1 \rangle $ to be $\deg_{1,2}(V_{0})=1$, $\deg_{1,2}(V_{1})=2$.
\end{defin}

With this definition we can obtain interesting results.

\begin{prop}\label{superindex}
	If $(V_0 , V_1)\in M_{N}(\complex)^2$ satisfies
	\begin{equation}\label{index}
	V_{1}^{i_{1}}V_{0}^{i_{2}}\cdots. V_{1}^{i_{n}}V_{0}^{i_{n+1}}=0,
	\end{equation}
	for any 	$i_{1}\geq 1$ ,
$	i_{1}+\cdots+i_{n+1}\leq n+1$, and $ n+2\leq \deg_{1,2}(V_{1}^{i_{1}}V_{0}^{i_{2}}\cdots. V_{1}^{i_{n}}V_{0}^{i_{n+1}})$,
 then 
	\begin{equation*}
	V_{j_{1}}\cdots V_{j_{k}}=0, \hspace{0.1cm}  j_{1}+\cdots+j_{k}\geq n, j_{1}\geq 1, k\geq 2.
	\end{equation*}

\end{prop}

\begin{demos}
Note that \eqref{index} implies that the monomials of degree $\geq n+1$ that begins with $V_{1}$ are zero. Since $j_{1}\geq 1$ we have that every monomial in $V_{j_{1}}\cdots V_{j_{k}}$ begins with $V_{1}$. Furthermore,  this polynomial is quasihomogenenous of degree $\deg_{1,2}(V_{j_{1}}\cdots V_{j_{k}})=(j_{1}+1)+ \cdots + (j_{k}+1)= j_{1}+\cdots +j_{k}+k\geq n+k \geq n+2$. In particular, the polynomial $V_{j_{1}}\cdots V_{j_{k}}$  is a linear combination of monomials of the form 
\eqref{index} which are zero. \QED
\end{demos}

\begin{remark}
	Two important elements satisfying \eqref{index} are $V_{1}^{i_{1}}$ and $V_{1}V_{0}^{n}$ with 
	\begin{equation*}
	i_{1} = \left \{
	\begin{array}{ll}
	\frac{n+2}{2}     & \mathrm{if\ } n \hspace{0.1cm} \text{is even} \\
\corch{\frac{n}{2}}+2 & \mathrm{if\ } n  \hspace{0.1cm} \text{is odd}\\
	\end{array}
	\right.
	\end{equation*}
\end{remark}

\begin{lemma}\label{subindex}
	For every $k\geq 0$, $((A_{1}^{[m]})^{k}P_{1}^{n+1}(V))_{i}$ is a polynomial in $V_{0}, V_{1}, \cdots ,V_{n}$ such that the sum of the subindices in its monomials is $\geq i$, $1\leq i \leq n+1$.
\end{lemma}
\begin{demos}
	The proof is by induction over $k$.
	For $k=0$ we are okay since
	\begin{equation*}
	(P_{1}^{n+1}(V))_{i}=P_{i}(V)=iV_{i}.
	\end{equation*}

Assume the claim for $k\geq 0$ and consider the case $k+1$ 
\begin{equation*}
((A_{1}^{[m]})^{k+1}P_{1}^{n+1}(V))_{i}=\sum_{j=1}^{n+1}(A_{1}^{[n]})_{ij}((A_{1}^{[m]})^{k}P_{1}^{n+1}(V))_{j}
\end{equation*}
\begin{equation*}
=\sum_{j=1}^{i}\frac{V_{i-j}}{2}((A_{1}^{[m]})^{k}P_{1}^{n+1}(V))_{j}+i((A_{1}^{[m]})^{k}P_{1}^{n+1}(V))_{i+1}.
\end{equation*}
Since the sum of the subindices of the monomials of $((A_{1}^{[m]})^{k}P_{1}^{n+1}(V))_{j}$ is $\geq j$ we obtain that the sum of the subindices in the monomials of 
$((A_{1}^{[m]})^{k+1}P_{1}^{n+1}(V))_{i}$ is $\geq \min\llave{(i-j)+j,i+1}=i$. \QED
\end{demos}

\begin{teor}[Bispectral Property  for a Class of  Polynomial Potentials]\label{potential}
	If $V(V_{0},V_{1},x)$ is a polynomial of degree $n$ such that $V''(x)=V'(x)V(x)$ and $(V_0 , V_1)\in M_{N}(\complex)^2$ satisfy \eqref{index}, then $V\in \mathbb{A}$. In particular, the operator $L=-\partial_{x}^2+V'(x)$ is bispectral.
\end{teor}
\begin{demos}
	Since $V''(x)=V'(x)V(x)$ we have that $L=-\partial_{x}^2+V'(x)$ satisfies $ (L\psi)(x,z)=-z^2 \psi(x,z)$ with 
	\begin{equation*}
	\psi (x,z)=\paren{Iz+\frac{1}{2}V(x)}e^{xz}.
	\end{equation*}

	On the other hand, from the Theorem \ref{bispectralth} we have $\mathbb{A}=\Gamma$. Furthermore, by the Lemma \ref{description} we have that the right bispectral algebra is the set of all $\theta\in M_{N}(\complex)[x]$  such that 
	$\corch{\theta,V}$ is a polynomial of degree $\leq m$ and 
	\begin{equation}\label{system}
	V_{-1}e_{1}(A_{1}^{[m]})^{k} P_{1}^{m+1}(\theta)=0, 
A_{2}^{[m]}(A_{1}^{[m]})^{k}P_{1}^{m+1}(\theta)=0, \hspace{0.1cm}	\text{for} \hspace{0.1cm} 0\leq k \leq (m+1)N-1,
	\end{equation}
	with $m=\deg(\theta)$ .

	However, since $V(x)$ is a polynomial  we have $V_{-1}=0$ and $\corch{V,V}=0$ is a polynomial of degree $\leq n:=\deg(V)$. 
	
	Note that 
	
	\begin{equation*}
(A_{2}^{[n]}(A_{1}^{[m]})^{k}P_{1}^{n+1}(V))_{i}
	=\sum_{j=1}^{n+1}(A_{2}^{[n]})_{ij}((A_{1}^{[m]})^{k}P_{1}^{n+1}(V))_{j}
	=\sum_{j=1}^{n+1}V_{i+n+1-j}((A_{1}^{[m]})^{k}P_{1}^{n+1}(V))_{j}.
	\end{equation*}

By the Lemma \ref{subindex} the sum of the subindices in the monomials of the polynomials $((A_{1}^{[m]})^{k}P_{1}^{n+1}(V))_{j}$ is $\geq j$. Therefore, the sum of the subindices of the monomials in the polynomial 
$(A_{2}^{[n]}(A_{1}^{[m]})^{k}P_{1}^{n+1}(V))_{i}$ is $\geq n$, $1\leq i \leq n+1$. Thus, the Proposition \ref{superindex} implies
\begin{equation*}
A_{2}^{[n]}(A_{1}^{[m]})^{k}P_{1}^{n+1}(V)=0, \hspace{0.1cm} k\in \natu.
\end{equation*}	
Then, $V\in \mathbb{A}$.	\QED
\end{demos}

We conclude this chapter with some examples applying the previous theorems.

\section{Illustrative Examples}

In this section we give  some ilustrative examples of bispectral operators $L=-\partial_{x}^2+V'(x)$ with polynomial potentials $V$ through the Theorem \ref{potential}.

\subsection{The Bispectral Algebra Associated to the Potential with Invertible Residue at \\ $x=0$ }

If $V_{-1}=-2I_{N}$, then $V_{j}=0$ for every $j\in \natu$ and $V(x)=\frac{-2I_{N}}{x}$. This implies that for every $m\in \natu$, 
\begin{equation*}
A_{1}^{[m]}=\paren{
	\begin{matrix} 
	0 & 0 & 0 &   \cdots & 0 & 0 & 0  \\
	0 & 0 & I_{N} &  \cdots & 0 & 0 & 0 \\
	. & . & \cdots & .  & \cdots & . & . \\
	. & . & \cdots & . & \cdots & . & . \\
	. & . & \cdots & . & \cdots & . & . \\
	0 &  0 & 0 &  \cdots & 0 & (m-2)I_{N}  & 0 \\
	0 & 0 & 0  &  \cdots &  0  & 0 & (m-1)I_{N} \\
	0 & 0 & 0  &  \cdots & 0 & 0 & 0\\
	\end{matrix}}
\end{equation*}
and $A_{2}^{[m]}=0$. 

Note that $A_{1}^{[m]}=\sum_{j=2}^{m}(j-1)I_{N}e_{j,j+1}$. We claim that $(A_{1}^{[m]})^{k}=\sum_{j=2}^{m-k+1}(j-1)j\cdot \cdot \cdot (j+k-2)I_{N}e_{j,j+k}$. We prove the claim by induction.
The initial step $k=1$ is clear. Assume $k\geq 1$ and note that 
\begin{equation*}
(A_{1}^{[m]})^{k+1}=(A_{1}^{[m]})^{k}A_{1}^{[m]}
=\paren{\sum_{j=2}^{m-k+1}(j-1)j\cdot \cdot \cdot (j-k+2)I_{N}e_{j,j+k}}\paren{\sum_{j=2}^{m}(j-1)I_{N}e_{j,j+1}}
\end{equation*}
\begin{equation*}
=\sum_{j=2}^{m-k+1}\sum_{l=2}^{m}(j-1)j\cdot \cdot \cdot (j+k-2)(l-1)I_{N}e_{j,j+k}e_{l,l+1}
\end{equation*}
\begin{equation*}
=\sum_{j=2}^{m-k}(j-1)j\cdot \cdot \cdot (j+k-2)(j+k-1)I_{N}e_{j,j+k}.
\end{equation*}

The claim follows by induction. We can write $(A_{1}^{[m]})^{k}=\sum_{j=2}^{m-k+1}\frac{(j+k-2)!}{(j-2)!}  I_{N}e_{j,j+k}$. In particular, 
$(A_{1}^{[m]})^{m-1}=(m-1)!e_{2,m+1}$ and $(A_{1}^{[m]})^{k}=0$ for every $k\geq m$.

This implies that, 
\begin{equation*}
(A_{1}^{[m]}+z)^{-1}=\sum_{k=0}^{\infty}\frac{(-1)^k (A_{1}^{[m]})^k}{z^{k+1}}
=\frac{I_{N}}{z}+\sum_{k=1}^{m-1}\sum_{j=2}^{m-k+1}(-1)^{k}\frac{(j+k-2)!I_{N}}{(j-2)!z^{k+1}}e_{j,j+k}.
\end{equation*}

Therefore, if $\theta(x)=\sum_{l=0}^{m}a_{l}x^l$ we have $P_{l}(\theta)=la_{l}$ and 
\begin{equation*}
c(z)=-(A_{1}^{[m]}+z)^{-1}P_{1}^{m+1}(\theta)=-\paren{\frac{I_{N}}{z}+\sum_{k=1}^{m-1}\sum_{j=2}^{m-k+1}(-1)^{k}\frac{(j+k-2)!I_{N}}{(j-2)!z^{k+1}}e_{j,j+k}}
\paren{\sum_{l=1}^{m+1}P_{l}(\theta)e_{l}}
\end{equation*}
\begin{equation*}
=\sum_{k=1}^{m-1}\sum_{l=k+2}^{m+1}\frac{(l-2)!(-1)^{k+1}}{(l-k-2)! z^{k+1}}P_{l}(\theta)e_{l-k}-\sum_{l=1}^{m+1}\frac{P_{l}(\theta)}{z}e_{l}
\end{equation*}
\begin{equation*}
=\sum_{k=1}^{m-1}\sum_{l=k+2}^{m+1}\frac{(l-2)!(-1)^{k+1}}{(l-k-2)! z^{k+1}}la_{l}e_{l-k}-\sum_{l=1}^{m+1}\frac{la_{l}}{z}e_{l}
\end{equation*}
\begin{equation*}
=\sum_{k=1}^{m-1}\sum_{s=2}^{m-k+1}\frac{(s+k-2)!(s+k)(-1)^{k+1}}{(s-2)! z^{k+1}}a_{s+k}e_{s}-\sum_{l=1}^{m+1}\frac{la_{l}}{z}e_{l}
\end{equation*}
\begin{equation*}
=\sum_{s=2}^{m}\paren{\sum_{k=1}^{m-s+1}\frac{(s+k-2)!(s+k)(-1)^{k+1}}{(s-2)! z^{k+1}}a_{s+k}}e_{s}-\sum_{l=1}^{m+1}\frac{la_{l}}{z}e_{l}
\end{equation*}
\begin{equation*}
=\sum_{s=2}^{m}\paren{\sum_{j=s+1}^{m}\frac{(j-2)!j(-1)^{j-s+1}}{(s-2)! z^{j-s+1}}a_{j}}e_{s}-\sum_{l=1}^{m+1}\frac{la_{l}}{z}e_{l}
\end{equation*}
\begin{equation*}
=-\frac{a_{1}}{z}e_{l}+\sum_{s=2}^{m}\paren{\sum_{j=s+1}^{m}\frac{(j-2)!j(-1)^{j-s+1}}{(s-2)! z^{j-s+1}}a_{j}-\frac{sa_{s}}{z}}e_{s}
\end{equation*}
\begin{equation*}
=-\frac{a_{1}}{z}e_{l}+\sum_{s=2}^{m}\paren{\sum_{j=s}^{m}\frac{(j-2)!j(-1)^{j-s+1}}{(s-2)! z^{j-s+1}}a_{j}}e_{s}=\sum_{s=1}^{m+1}c_{s-1}(z)e_{s}.
\end{equation*}

Thus, $c_{0}(z)=b_{0}(z)-a_{0}=-\frac{a_{1}}{z}$, $b_{0}(z)=a_{0}-\frac{a_{1}}{z}$. Furthermore, 
\begin{equation*}
c_{s-1}(z)=\sum_{j=s}^{m}\frac{(j-2)!j(-1)^{j-s+1}}{(s-2)! z^{j-s+1}}a_{j}, 2\leq s \leq m+1.
\end{equation*}
In other words, 
\begin{equation*}
c_{s}(z)=\sum_{j=s+1}^{m}\frac{(j-2)!j(-1)^{j-s}}{(s-1)! z^{j-s}}a_{j}, 1\leq s \leq m.
\end{equation*}
On the other hand, we have the restrictions $\frac{1}{2}V_{-1}c_{0}(z)=0$ and $P_{0}(\theta)=0$. In this case the former restriction says that $c_{0}(z)=0$ and the last is redundant. Therefore, $b_{0}(z)=a_{0}$ and $a_{1}=0$.

We conclude that, $\mathbb{A}=\llave{\theta\in M_{N}(\complex[x])\mid \theta^{'}(0)=0}$ and for every $\theta\in \mathbb{A}$,
$\theta(x)=\sum_{l=0}^{m}a_{l}x^l$  there exists $B(z,\partial_{z})=
\sum_{k=0}^{m}\partial_{z}^{k} \cdot b_{k}(z)= a_{0}+\partial_{z} \cdot \sum_{j=2}^{m}\frac{(j-2)!j(-1)^{j-1}}{z^{j-1}}a_{j}+\sum_{k=2}^{m}\partial_{z}^{k} \cdot \llave{a_{k}+\sum_{j=k+1}^{m}\frac{(j-2)!j(-1)^{j-k}}{(k-1)! z^{j-k}}a_{j} }$ such that 
$(\psi B)(x,z)=\theta(x)\psi(x,z)$.

\subsection{Examples of Polynomial Potentials of  degree $n=1,2,3$}

In this subsection we use Maxima to perform explicit examples of Theorem \ref{potential}.

\begin{itemize}
	\item For $n=1$ the equations \eqref{index} turns out to be $V_{1}V_{0}=V_{1}^2=0$. 
	For $N=2$ we consider
	\begin{equation*}
	V_{1}=\paren{
		\begin{matrix} 
		0 & 1 \\
		0 & 0 \\
		\end{matrix}}
	\end{equation*}
	and 
	\begin{equation*}
	V_{0}=\paren{
		\begin{matrix} 
		V_{011} & V_{012} \\
		0 & 0 \\
		\end{matrix}}.
	\end{equation*}
	To obtain the potential $V(x)=V_0 + V_1 x$.
	\item  For $n=2$ the equations \eqref{index} turns out to be $V_{1}V_{0}V_{1}=V_{1}V_{0}^2=V_{1}^2=0$.
	For $N=4$ we consider
	\begin{equation*}
	V_{0}=\paren{
		\begin{matrix} 
		V_{011} & V_{012} & V_{013} & V_{014}  \\
		0 & 0 & V_{023} & V_{024}\\
			0 & 0 & 0 & 0\\
					0 & 0 & 0 & 0\\
		\end{matrix}}
	\end{equation*}
	and 
		\begin{equation*}
	V_{1}=\paren{
		\begin{matrix} 
		0 & V_{112} & 0 & 0  \\
		0 & 0 & 0 & 0\\
		0 & 0 & 0 & 0\\
		0 & 0 & 0 & 0\\
		\end{matrix}}.
	\end{equation*}
	To obtain 
		\begin{equation*}
	V_{2}=\paren{
		\begin{matrix} 
		0 & 0 & \frac{V_{023}V_{112}}{2} & \frac{V_{024}V_{112}}{2}  \\
		0 & 0 & 0 & 0\\
		0 & 0 & 0 & 0\\
		0 & 0 & 0 & 0\\
		\end{matrix}}
	\end{equation*}
	and the potential $V(x)=V_{0}+V_{1}x+V_{2}x^2$.
	\item  For $n=3$ the equations \eqref{index} turns out to be $ V_{1}^3=V_{1}V_{0}V_{1}=V_{1}V_{0}^3 =V_{1}^2V_{0}=V_{1}V_{0}^2 V_{1}$.
	For $N=4$ we consider
		\begin{equation*}
	V_{0}=\paren{
		\begin{matrix} 
		V_{011} & V_{012} & V_{013} & V_{014}  \\
		0 & 0 & V_{023} & V_{024}\\
		0 & 0 & 0 & 0\\
		0 & 0 & 0 & 0\\
		\end{matrix}}
	\end{equation*}
	and 
	\begin{equation*}
	V_{1}=\paren{
		\begin{matrix} 
			V_{111} & V_{112} & V_{113} & V_{114}  \\
		0 & 0 & V_{123} & V_{124}\\
		0 & 0 & 0 & 0\\
		0 & 0 & 0 & 0\\
		\end{matrix}}.
	\end{equation*}
	To obtain 
		\begin{equation*}
	V_{2}=\paren{
		\begin{matrix} 
		0 & 0 & \frac{V_{023}V_{112}}{2} & \frac{V_{024}V_{112}}{2}  \\
		0 & 0 & 0 & 0\\
		0 & 0 & 0 & 0\\
		0 & 0 & 0 & 0\\
		\end{matrix}}
	\end{equation*}
	and 
		\begin{equation*}
	V_{3}=\paren{
		\begin{matrix} 
		0 & 0 & \frac{V_{112}V_{123}}{2} & \frac{V_{112}V_{124}}{2}  \\
		0 & 0 & 0 & 0\\
		0 & 0 & 0 & 0\\
		0 & 0 & 0 & 0\\
		\end{matrix}}
	\end{equation*}
		and the potential $V(x)=V_{0}+V_{1}x+V_{2}x^2+V_{3}x^3$.
	
\end{itemize}

\section*{Conclusions and Comments}

In this article, we characterized the bispectral algebra associated with some type of matrix Schr\"{o}dinger operators with polynomial potential. This characterization was achieved using the family of functions 
$\mathcal{P}=\llave{P_{k}}_{k\in \natu}$. The natural question arising in the characterization of the bispectral algebra associated with a matrix Schr\"{o}dinger operator is about the possibility of writing its conditions using some family of maps  $\mathcal{P}$ satisfying the Leibniz rule.

\subsubsection*{Acknowledgments}

BDVC acknowledges the support from the FSU2022-010 grant from Khalifa University, UAE.

\section*{ORCID}

 Brian Vasquez Campos's ORCID: {0000-0003-3922-9956}

	\end{document}